\documentclass{amsart}

\title{Congruences of elliptic curves arising from\\ non-surjective mod $N$ Galois representations}

\author{Sam Frengley} 
\address{University of Cambridge, DPMMS, Centre for Mathematical Sciences, Wilberforce Road, Cambridge CB3 0WB, UK}
\email{stf32@cam.ac.uk}

\date{June 16, 2022}

\usepackage{graphicx}

\usepackage{fullpage}	
\usepackage{colonequals}
\usepackage[pagebackref=true,hypertexnames=false]{hyperref}

\usepackage{amsmath}	
\usepackage{amssymb}	
\usepackage{amsthm}		
\usepackage{enumitem}	

\usepackage{theoremref}

\usepackage{tikz-cd}

\usepackage{multirow}
\usepackage{diagbox}
\usepackage{adjustbox}
\usepackage{float}
\usepackage{arydshln}
\usepackage{makecell}
\usepackage{pdflscape}

\numberwithin{table}{section}

\newcommand{\bbF}{\mathbb{F}}
\newcommand{\bbP}{\mathbb{P}}
\newcommand{\bbQ}{\mathbb{Q}}
\newcommand{\bbZ}{\mathbb{Z}}

\DeclareMathOperator{\Gal}{Gal}
\DeclareMathOperator{\Aut}{Aut}
\DeclareMathOperator{\Tr}{tr}

\DeclareMathOperator{\GL}{GL}

\DeclareMathOperator{\PGL}{PGL}
\DeclareMathOperator{\ns}{ns}
\DeclareMathOperator{\s}{s}

\newcommand{\LMFDBLabel}[1]{\textnormal{\href{https://www.lmfdb.org/EllipticCurve/Q/#1/}{\texttt{#1}}}}
\newcommand{\LMFDBLabelIsog}[1]{\textnormal{\href{https://www.lmfdb.org/EllipticCurve/Q/#1/}{\texttt{#1*}}}}

\newcommand{\eniii}{\textnormal{(\roman*)}}
\newcommand{\enABC}{\textnormal{(\Alph*)}}

\newtheorem{theorem}{Theorem}
\numberwithin{theorem}{section} 
\newtheorem{prop}[theorem]{Proposition}
\newtheorem{coro}[theorem]{Corollary}  
\newtheorem{lemma}[theorem]{Lemma}
\newtheorem{conj}[theorem]{Conjecture}

\theoremstyle{definition}\newtheorem{defn}[theorem]{Definition}
\theoremstyle{definition}\newtheorem{remark}[theorem]{Remark}
\theoremstyle{definition}\newtheorem{example}[theorem]{Example}

\begin{document}
	\begin{abstract}
		We study $N$-congruences between quadratic twists of elliptic curves. If $N$ has exactly two distinct prime factors we show that these are parametrised by double covers of certain modular curves. In many, but not all cases, the modular curves in question correspond to the normaliser of a Cartan subgroup of $\GL_2(\bbZ/N\bbZ)$. By computing explicit models for these double covers we find all pairs $(N, r)$ such that there exist infinitely many $j$-invariants of elliptic curves $E/\bbQ$ which are $N$-congruent with power $r$ to a quadratic twist of $E$. We also find an example of a $48$-congruence over $\bbQ$. We make a conjecture classifying nontrivial $(N,r)$-congruences between quadratic twists of elliptic curves over $\bbQ$.
		
		Finally, we give a more detailed analysis of the level $15$ case. We use elliptic Chabauty to determine the rational points on a modular curve of genus $2$ whose Jacobian has rank $2$ and which arises as a double cover of the modular curve $X(\ns 3^+, \ns 5^+)$. As a consequence we obtain a new proof of the class number $1$ problem.
	\end{abstract}

	\maketitle
	
	\section{Introduction}
	
		Let $K$ be a field of characteristic $0$ and let $G_K \colonequals \Gal(\bar{K}/K)$ be the absolute Galois group of $K$. Let $E$ be an elliptic curve defined over $K$ and for each integer $N \geq 1$, let $E[N]$ denote the $N$-torsion subgroup of $E(\bar{K})$.
		
		We say that a pair of elliptic curves $E/K$ and $E'/K$ are \emph{$N$-congruent} (over $K$) if there exists an isomorphism $\phi \colon E[N] \to E'[N]$ of $G_K$-modules. We say that an $N$-congruence $\phi$ has \emph{power $r$} if $\phi$ raises the Weil pairing to a power of $r$. By composing $\phi$ with the multiplication by $k$ map on $E'[N]$ we obtain an $N$-congruence with power $k^2 r$, hence $r$ is defined up to multiplication by a square in $(\bbZ/N\bbZ)^\times$. For brevity we will say that an $N$-congruence with power $r$ is an \emph{$(N, r)$-congruence}. 
		
		If $E$ and $E'$ are $m$-isogenous over $K$, and $(N,m) =1$, then the restriction of the isogeny to the $N$-torsion is an $(N,m)$-congruence. We say that a congruence is \emph{trivial} if $E$ and $E'$ are isogenous over $K$. Note that, for simplicity, we do not require a trivial $N$-congruence to be equal to the restriction of the isogeny to $E[N]$ (for further discussion see \thref{remark:weirdtrivialcongs}). 

		The study of congruences has been motivated by, for example, Diophantine problems (in particular, solutions to generalised Fermat equations, see e.g., \cite{PSS_TOX7APSTX2Y3Z7}) and the study of genus $2$ curves with split Jacobians (see e.g., \cite{BD_TAOG2CW44SJ}).

		For small values of $N$, success has been had finding examples of $N$-congruences by computing twists of the modular curve $X(N)$ which parametrise elliptic curves $(N, r)$-congruent to a fixed elliptic curve $E$. For $N \geq 7$ the curve $X(N)$ has genus $\geq 2$ and it has been more fruitful to compute so-called ``modular diagonal quotient surfaces'' whose points parametrise pairs of $(N, r)$-congruent elliptic curves. When $K = \bbQ$, these approaches have led to a proof that for each $N \leq 13$ there exist infinitely many pairs of nontrivially $N$-congruent elliptic curves (see \cite{F_OFO13CEC} for the case $N = 13$, and \cite{F_EMSFCOEC} and the references therein for the case $N \leq 12$ - note that not all powers have been considered when $N = 12$). 
		
		On the other hand, it is a conjecture of Frey and Mazur that, over $\bbQ$, nontrivial $N$-congruences between elliptic curves should not exist for sufficiently large $N$. When $N = p$ is prime, this has been refined by Fisher who has conjectured that there exist no pairs of nontrivially $p$-congruent elliptic curves for $p > 17$ (see \cite[Conjecture 1.1]{F_OPO17CEC}). Cremona and Freitas \cite[Theorem 1.3]{CF_GMFTSTOCBEC} have determined all pairs of nontrivially $p$-congruent elliptic curves with conductor less than $500,000$. Specifically they showed that when $p > 17$ there exist no pairs of non-isogenous, $p$-congruent, elliptic curves with conductor $<500,000$. 

		It was observed by Halberstadt \cite[Section 2]{H_SLCMX11} that an elliptic curve such that the image of the corresponding mod $p$ Galois representation, $\bar{\rho}_{E, p}$, is contained in the normaliser of a Cartan subgroup of $\GL_2(\bbF_p)$, but not the Cartan subgroup itself, admits a $p$-congruence with a nontrivial quadratic twist. An elliptic curve is $2$-congruent to every quadratic twist, hence such a congruence is automatically a $2p$-congruence. Based upon this observation Halberstadt \cite{H_SLCMX11} found examples of nontrivial $10, 14,$ and $22$-congruences of elliptic curves over $\bbQ$ by using explicit models for $X(\ns  p^+)$ for $p = 5,7,11$ (indeed in each of these cases there are infinitely many such congruences, see \thref{thm:infmany}).

		\begin{defn}
			We say that an $(N, r)$-congruence between elliptic curves $E/K$ and $E'/K$ is of \emph{quadratic twist-type} if $E$ is a nontrivial quadratic twist of $E'$. In this case we say that $E$ (and $E'$) \emph{admits a quadratic twist-type $(N,r)$-congruence}.

			We say that a quadratic twist-type congruence is \emph{nontrivial} if the congruence itself is nontrivial (i.e, if $E$ and $E'$ are not isogenous over $K$). 
		\end{defn}

		Cremona--Freitas \cite[Theorem 2.4]{CF_GMFTSTOCBEC} proved that (under their condition ``\textbf{(S)}'', which is satisfied when $K = \bbQ$ and $p \geq 7$) quadratic twist-type $p$-congruences occur \emph{if and only if} the image of the corresponding mod $p$ Galois representation is contained in the normaliser of a Cartan subgroup of $\GL_2(\bbF_p)$, but not the Cartan subgroup itself. 
		
		By a different proof we extend this result to odd prime powers (see \thref{cartan_cong_pl}). Our result is that if $p$ is an odd prime and $l \geq 1$ is an integer, then an elliptic curve $E/K$ admits a quadratic twist-type $p^l$-congruence if and only if the image of corresponding mod $p^l$ Galois representation is contained in the normaliser of a Cartan subgroup of $\GL_2(\bbZ/p^l\bbZ)$, but not the Cartan subgroup itself. Note that we avoid condition ``\textbf{(S)}''.
		
		Therefore classifying quadratic twist-type $p^l$-congruences (hence $2p^l$-congruences) is equivalent to finding all (non-cuspidal) rational points on the modular curves $X(\s {p^l}^+)$ and $X(\ns {p^l}^+)$ which do not lift to points on $X(\s {p^l})$ respectively $X(\ns {p^l})$. 
		
		If $K = \bbQ$ and $p^l > 7$ then CM points correspond to trivial congruences (see \thref{lemma:CM_congs} and \thref{prop:nontrivialCMcongs}). Therefore classifying nontrivial quadratic twist-type $p^l$-congruences is equivalent to finding all \emph{exceptional} (i.e., non-cuspidal, non-CM) rational points on the modular curves $X(\s {p^l}^+)$ and $X(\ns {p^l}^+)$ which do not lift to points on $X(\s {p^l})$ respectively $X(\ns {p^l})$. The curves $X(\s {p^l}^+)$ have no exceptional points for $p^l > 7$ (see \cite[Theorem 1.1]{BPR_RPOX+0pr} for $p > 7$ except $p = 13$, \cite{BDMSTV_ECKFTSCMCOL13} in the level $13$ case, and \cite{RSZB_LAIOGFECOQ} when $ p^l = 9, 25, 49$). It is a well known open problem to determine the rational points on $X(\ns {p^l}^+)$ when $p^l \geq 19$ (when $p^l = 13$ and $17$ this has been resolved in \cite{BDMSTV_ECKFTSCMCOL13} and \cite{BDMTV_QCFMCAAE} respectively).

		We extend the results of Halberstadt and Cremona--Freitas to composite $N$, and in particular we prove:

		\begin{theorem} \thlabel{thm:infmany}
			There are infinitely many $j$-invariants of elliptic curves $E/\bbQ$ admitting nontrivial quadratic twist-type $(N,r)$-congruences if and only if $(N,r) \in \mathcal{S}$ where
			\begin{equation*}
				\mathcal{S} = \left\{ 
				\begin{aligned}
					&(2,1), (3,1), (3,2), (4,1), (4,3), (5,1), (5,2), (6,1), (6,5), (7,1), (7,3), 
					\\ &(8,1), (8,3), (8,5), (8,7), (9,1), (10,1), (10,3), (11,1), (12,7), (12,11), 
					\\ & (14,1), (14,3), (16,3), (18,1), (20,3), (22,1), (24,19), (28,3), (28,11), (36,19)
				\end{aligned}
				\right\}
			\end{equation*}
			and the second entry is viewed as an element of $(\bbZ/N\bbZ)^\times$ modulo squares.
		\end{theorem}

		In the spirit of Mazur's ``Program B'' (and other work on modular curves - see e.g., \cite{RZB_ECOQA2AIOG}, \cite{SZ_MCOPPLWIMRP}, and \cite{RSZB_LAIOGFECOQ}) it would be interesting to extend \thref{thm:infmany} to classify \emph{all} nontrivial quadratic twist-type $(N, r)$-congruences (assuming the conjecture that the modular curves $X(\ns {p^l}^+)$ have no exceptional points for $p^l \geq 19$, see e.g., \cite[Conjecture 1.5]{RSZB_LAIOGFECOQ}).

		Towards this goal we exhibit the following examples of nontrivial quadratic twist-type congruences not arising from the families in \thref{thm:infmany}.  

		\begin{theorem} \thlabel{thm:examples}
			Let $E_{(30,7)}$, $E_{(32,3)a}$, $E_{(32,3)b}$, and $E_{(48,19)}$ be the elliptic curves defined over $\bbQ$ by the Weierstrass equations
			\begin{align*}
				E_{(30,7)} &: y^2 + xy = x^3 - x^2 - 273176601587417x - 1741818799948905109620,\\
				E_{(32,3)a} &: y^2 + y = x^3 - 3014658660x + 150916472601529, \\
				E_{(32,3)b} &:  y^2 + y = x^3 - 26668521183591560x - 1676372876256754456631251, \\
				E_{(48,19)} &: y^2 + y = x^3 + 468240736152891010x - 148374586624464876247316957. 
			\end{align*}
			Then $E_{(30,7)}$ is $(30,7)$-congruent (over $\bbQ$) to its quadratic twist by $-214663$, and $E_{(32,3)a}$, $E_{(32,3)b}$, and $E_{(48,19)}$ are $(32,3)$, $(32,3)$, and $(48,19)$-congruent respectively to their quadratic twists by their discriminants.
		\end{theorem}

		The mod $N$ Galois representations of a pair of $N$-congruent elliptic curves $E/\bbQ$ and $E'/\bbQ$ are isomorphic, hence have the same trace. Therefore for any good prime $\ell$ for both $E$ and $E'$ which does not divide $N$ we have $a_{\ell}(E) \equiv a_{\ell}(E') \pmod{N}$. By computing traces of Frobenius we see that none of the $N$-congruences in \thref{thm:examples} is an $M$-congruence for any $M > N$.

		By searching for points on the relevant modular curves of small genus we give evidence for the following conjecture:

		\begin{conj} \thlabel{conj:allqtt}
			If $E/\bbQ$ is an elliptic curve admitting a nontrivial quadratic twist-type $(N,r)$-congruence then either $(N,r) \in \mathcal{S}$ or $E$ is a quadratic twist of either $E_{(30,7)}$, $E_{(32,3)a}$, $E_{(32,3)b}$, or $E_{(48,19)}$.
		\end{conj}

		The paper is arranged as follows:

		We first give a short description of Cartan subgroups of $\GL_2(\bbZ/N\bbZ)$ and modular curves in Section \ref{sec:basic_defs}.

		In Section \ref{sec:technical} we prove our main technical results. By adapting the proofs of Cremona--Freitas \cite[Theorem 2.4]{CF_GMFTSTOCBEC} and Halberstadt \cite[Section 2.1]{H_SLCMX11} we give a necessary and sufficient condition for an elliptic curve to admit a quadratic twist-type $N$-congruence for $N > 2$, in terms of the mod $N$ image of Galois (see \thref{cong_q_twist}). We first use this to extend \cite[Theorem 2.4]{CF_GMFTSTOCBEC} to odd prime powers (see \thref{cartan_cong_pl}). We then classify quadratic twist-type $(2^k, r)$-congruences of non-CM elliptic curves as an application of Rouse and Zureick-Brown's \cite{RZB_ECOQA2AIOG} classification of the possible $2$-adic images of the Galois representation attached to an elliptic curve $E/\bbQ$ without CM. 

		In Section \ref{sec:trivialCM} we show that an elliptic curve $E/\bbQ$ is isogenous to a nontrivial quadratic twist if and only if $E$ has CM (see \thref{lemma:CM_congs}). By classifying all nontrivial quadratic twist-type congruences between elliptic curves, $E/\bbQ$, with CM (see \thref{prop:nontrivialCMcongs}) we deduce \thref{thm:infmany} when $N$ is a prime power.

		In Section \ref{sec:doublecovers} we construct modular curves which parametrise elliptic curves admitting quadratic twist-type $N$-congruences for composite $N$. By computing the genera of these modular curves and computing explicit models for those with genus $\leq 1$ we conclude the proof of \thref{thm:infmany} in Section \ref{sec:2kpl}.

		In Section \ref{sec:bigconj} we give evidence for \thref{conj:allqtt}. Specifically we compute explicit models for modular curves of genus $\leq 4$ whose points give rise to quadratic twist-type congruences. By searching for rational points on these modular curves we prove \thref{thm:examples} and a conditional form of \thref{conj:allqtt} (see \thref{thm:conditional}).
		
		Finally in Section \ref{sec:15-7ratpts} we perform a detailed analysis of quadratic twist-type $(15,7)$-congruences (equivalently $(30,7)$-congruences). We use elliptic Chabauty to find all the rational points on a genus $2$ curve whose Jacobian has rank $2$ and admits a Richelot isogeny. It follows that $E_{(30,7)}/\bbQ$ is (up to a quadratic twist) the only elliptic curve over $\bbQ$ admitting a nontrivial quadratic twist-type $(15,7)$-congruence (see Section \ref{sec:ratpts}). As a corollary we obtain a new proof of the class number $1$ problem in Section \ref{sec:clnum}.		

		Our interest in this problem arose from extending the search of Cremona--Freitas for $p$-congruences of elliptic curves within the LMFDB (see \cite[Section 3]{CF_GMFTSTOCBEC}) to $N$-congruences for composite $N$. We observed that (except for some cases when $N = 16$ which are still to be resolved) when $N \geq 15$ all $N$-congruences between pairs of non-isogenous elliptic curves appearing within the LMFDB are of quadratic twist-type, with a single exception when $N = 17$ (this is the $17$-congruence in \cite[Section 3.7]{CF_GMFTSTOCBEC}).
		
		All computations were performed in \texttt{Magma} \cite{MAGMA}. We make extensive use of Rouse, Sutherland, and Zureick-Brown's package \texttt{gl2.m} (see \cite{RSZB_Electronic}) which is based on the techniques in \cite{RSZB_LAIOGFECOQ}. In particular we found this package especially useful for computing genera of modular curves and manipulating subgroups of $\GL_2(\bbZ/N\bbZ)$. The supporting code may be found in our GitHub repository \cite{ME_ELECTRONIC}.
		
		\subsection{Acknowledgements}
			 I would like to express my deep gratitude to Tom Fisher for introducing me to the topic and for many patient and enlightening discussions. I would also like to thank the Woolf Fisher and Cambridge Trusts for their financial support and an anonymous referee for several helpful comments on an earlier version of this paper.

	\section{Cartan Subgroups of \texorpdfstring{$\GL_2(\bbZ/p^l\bbZ)$}{GL\_2(Z/p{\textasciicircum}lZ)} and Modular Curves} \label{sec:basic_defs}
		Let $p$ be a prime number, and let $l \geq 1$ be an integer. We briefly summarise some facts about Cartan subgroups of $\GL_2(\bbZ/p^l\bbZ)$. Our exposition closely follows that of Baran \cite{B_NONSCSMCATCNP}.

		Let $\mathcal{O}$ be an order in a quadratic number field such that $p$ is unramified in $\mathcal{O}$. Consider the free, rank $2$, $\bbZ/p^l\bbZ$-algebra $A = \mathcal{O}/p^l\mathcal{O}$. Now $A/pA$ is either isomorphic to $\bbF_p \times \bbF_p$ or $\bbF_{p^2}$, in the first case we say that $A$ is \emph{split}, and in the second case we say $A$ is \emph{nonsplit}. Up to isomorphism there is a unique split (respectively nonsplit) $\bbZ/p^l\bbZ$-algebra $A$. 

		Let $A^\times$ act on $A$ by multiplication. By choosing a $\bbZ/p^l\bbZ$-basis for $A$ we obtain an embedding $A^\times \hookrightarrow \GL_2(\bbZ/p^l\bbZ)$. We define a \emph{Cartan subgroup} of $\GL_2(\bbZ/p^l\bbZ)$ to be any subgroup occurring as the image of such an embedding, and we say that a Cartan subgroup is \emph{split} (respectively \emph{nonsplit}) if $A$ is split (respectively nonsplit).
		
		The image of the representation $A^\times \hookrightarrow \GL_2(\bbZ/p^l\bbZ)$ depends only on the basis we choose for $A$. In particular all split (respectively nonsplit) Cartan subgroups of $\GL_2(\bbZ/p^l\bbZ)$ are conjugate in $\GL_2(\bbZ/p^l\bbZ)$. 

		We now describe how to explicitly construct Cartan subgroups. Explicit descriptions of these (and other) subgroups of $\GL_2(\bbZ/N\bbZ)$ may also be found in \cite[Section 3]{DGM_AOCMCOPACL}. 

		\subsection*{A Split Cartan Subgroup}

		Let $G(\s p^l) \subset \GL_2(\bbZ/p^l\bbZ)$ denote the subgroup of diagonal matrices. Then $G(\s p^l)$ is a split Cartan subgroup of $\GL_2(\bbZ/p^l\bbZ)$ since it is the regular representation of $A = (\bbZ/p^l\bbZ)^2$ with respect to the standard basis.

		When $p$ is odd, the normaliser $G(\s  {p^l}^+)$ of $G(\s p^l)$ in $\GL_2(\bbZ/p^l\bbZ)$ is generated by $G(\s p^l)$ and $\big(\begin{smallmatrix} 0 & 1\\ 1 & 0 \end{smallmatrix}\big)$. By convention for each integer $k \geq 1$ we write $G(\s  {2^k}^+)$ for the subgroup of $\GL_2(\bbZ/2^k\bbZ)$ generated by $G(\s 2^k)$ and $\big(\begin{smallmatrix} 0 & 1\\ 1 & 0 \end{smallmatrix}\big)$, note however that $G(\s  {2^k}^+)$ is \emph{strictly} contained in the normaliser of $G(\s 2^k)$ in $\GL_2(\bbZ/2^k\bbZ)$ since, for example, $\big(\begin{smallmatrix} 1 & 2^{k-1}\\ 0 & 1 \end{smallmatrix}\big)$ normalises $G(\s 2^k)$\footnote{This observation means that is often convenient to define the normaliser of a Cartan subgroup of $\GL_2(\bbZ/2^k\bbZ)$ to be a subgroup which contains a Cartan subgroup with index $2$ (see e.g., \cite[\href{https://www.lmfdb.org/knowledge/show/gl2.normalizer_cartan}{Normalizer of a Cartan subgroup}]{lmfdb}).}.
		
		Note that for each $p$ the subgroups $G(\s p^l)$ and $G(\s {p^l}^+)$ are equal to their transposes.

		\subsection*{A Nonsplit Cartan Subgroup}

		Let $p$ be prime and let $\mathcal{O} = \mathbb{Z}[\alpha]$ be a quadratic order such that $p$ is inert in $\mathcal{O}$. Let $f(x) = x^2 - ux + v \in \mathbb{Z}[x]$ be the minimal polynomial of $\alpha$. For the purposes of this article we take $f(x) = x^2 + x + 1$ when $p = 2$, and when $p$ is odd, we take $f(x) = x^2 - \xi$ where $\xi \in \bbZ$ is a quadratic nonresidue modulo $p$. 
		
		The $\bbZ/p^l\bbZ$-algebra $A = \mathcal{O}/p^l\mathcal{O}$ is nonsplit. Fix the $\bbZ/p^l\bbZ$-basis $\{1, \alpha\}$ for $A$, then $A^\times$ embeds in $\GL_2(\bbZ/p^l\bbZ)$ by its regular action on this basis. By construction the image is a nonsplit Cartan subgroup, which we denote $G(\ns  p^l)$. Explicitly for each integer $k \geq 1$ we take
		\begin{equation*}
			G(\ns  2^k) = \left\{ \begin{pmatrix} a & -b \\ b & a-b \end{pmatrix} : a,b \in \bbZ/2^k\bbZ \text{ such that $(a,b) \not\equiv (0,0) \pmod{2}$} \right\}
		\end{equation*}
		and when $p$ is odd we take
		\begin{equation*}
			G(\ns  p^l) = \left\{ \begin{pmatrix} a & b\xi \\ b & a \end{pmatrix} : a,b \in \bbZ/p^l\bbZ \text{ such that $(a,b) \not\equiv (0,0) \pmod{p}$} \right\}.
		\end{equation*}

		For any prime number $p$, the normaliser $G(\ns {p^l}^+)$ of $G(\ns  p^l)$ in $\GL_2(\bbZ/p^l\bbZ)$ is generated by $G(\ns  p^l)$ and the matrix $\big(\begin{smallmatrix} 1 & u\\ 0 & -1 \end{smallmatrix}\big)$ (i.e., $\big(\begin{smallmatrix} 1 & -1\\ 0 & -1 \end{smallmatrix}\big)$ when $p = 2$ and $\big(\begin{smallmatrix} 1 & 0\\ 0 & -1 \end{smallmatrix}\big)$ when $p$ is odd).
		
		The subgroups $G(\ns 2^k)$ and $G(\ns {2^k}^+)$ are equal to their transposes. When $p$ is odd the transpose of $G(\ns p^l)$ is again a nonsplit Cartan subgroup. In particular $G(\ns p^l)$ is conjugate to its transpose and the same is true for $G(\ns {p^l}^+)$.
		
		\subsection{Modular Curves}
			We give a brief account of the modular curve $X(H)$. For more details see e.g., \cite[Section 3]{B_NONSCSMCATCNP}, \cite[Section IV-3]{DR_LSDMDCE}, \cite[Section 2.3]{RSZB_LAIOGFECOQ}, \cite[Section 2]{RZB_ECOQA2AIOG}, \cite[A.5]{S_LOTMWT}, or \cite[Section 2]{SZ_MCOPPLWIMRP}.

			Let $N \geq 1$ be an integer, $K$ be a field of characteristic $0$, and $E/K$ be an elliptic curve. By fixing a $\bbZ/N\bbZ$-basis for $E[N]$, we identify $\Aut(E[N]) \cong \GL_2(\bbZ/N\bbZ)$ acting by matrix multiplication on the left. Let 
			\begin{equation*}
				\bar{\rho}_{E, N} \colon G_{K} \to \Aut(E[N]) \cong \GL_2(\bbZ/N\bbZ)
			\end{equation*}
			be the mod $N$ representation of the absolute Galois group $G_K$.
			
			Let $Y(N)/\bbQ$ denote the moduli space of elliptic curves with full level $N$-structure. That is, if $S/\bbQ$ is a scheme, the $S$-points on $Y(N)$ parametrise (equivalence classes of) pairs $(E/S, \iota)$, where $E/S$ is an elliptic curve, and $\iota \colon (\bbZ/N\bbZ)^2_S \to E[N]$ is an isomorphism of $S$-group schemes (i.e., $\iota$ chooses a $\bbZ/N\bbZ$-basis for $E[N]$). 
			
			The natural (left) action of a matrix $g \in \GL_2(\bbZ/N\bbZ)$ on $(\bbZ/N\bbZ)^2_S$ induces an automorphism of $Y(N)$ via $g \cdot (E/S, \iota) = (E/S, \iota \circ g)$. Denote by $X(N)/\bbQ$ the compactification of $Y(N)$. The action of $\GL_2(\bbZ/N\bbZ)$ on $Y(N)$ extends uniquely to an action on $X(N)$.

			\begin{remark}
				With this convention, $X(N)$ is \emph{not} geometrically irreducible, it consists of $\varphi(N)$ Galois conjugate components, permuted by $\Gal(\bbQ(\zeta_N)/\bbQ)$.
			\end{remark}
			
			Let $H \subset \GL_2(\bbZ/N\bbZ)$ be a subgroup and define $X(H)$ to be the quotient $X(N)/H$. For simplicity we will also assume that $-I \in H$ so that the quotient is a scheme (see the discussion after Remark 2.2 in \cite{RZB_ECOQA2AIOG}). Note that $X(H)$ is geometrically irreducible if and only if $\det(H) = (\bbZ/N\bbZ)^\times$. By forgetting the level structure we have a natural map $X(H) \to X(1) \cong \bbP^1$ defined over $\mathbb{Q}$ and ramified only above $j = 0, 1728, \infty$. We say that a point $P \in X(H)$ is a \emph{cusp} if it is in the fibre above $j = \infty$.

			A rational point $P$ on a modular curve $X(H)$ is said to be \emph{exceptional} if $X(H)$ has finitely many rational points and $P$ is neither a cusp nor lies in the fibre of $X(H) \to X(1)$ above a CM $j$-invariant (e.g., \cite[Definition 1.2]{RSZB_LAIOGFECOQ}). We say that a $j$-invariant is exceptional if it is the $j$-invariant of an exceptional point (i.e., the image of an exceptional point under the $j$-map $X(H) \to X(1)$). 

			\begin{remark} \thlabel{remark:conjugatetotranspose}
				Since we only consider modular curves $X(H)$ where the subgroup $H \subset \GL_2(\bbZ/N\bbZ)$ is conjugate to its own transpose, our choice of left (versus right) action is not important. For further discussion see \cite[Remark 2.2]{RZB_ECOQA2AIOG}.
			\end{remark}
			
			The purpose of this construction is the following lemma:
			\begin{lemma} \thlabel{lemma:modcurveimages}
				Let $K$ be a field of characteristic $0$ and let $E/K$ be an elliptic curve such that $j(E) \neq 0, 1728$. Then there exists a $K$-rational point in the fibre of $X(H) \to X(1)$ above $j(E)$ if and only if the image, $\bar{\rho}_{E, N}(G_K)$, of the corresponding mod $N$ Galois representation is contained in a subgroup conjugate to $H$.

				\begin{proof}
					See \cite[Lemma 2.1]{RZB_ECOQA2AIOG} and \cite[Proposition 3.2]{Z_PIOTGIOECOQ}.
				\end{proof}
			\end{lemma}

			From the properties of the Weil pairing it follows that if $H$ occurs as the image of the mod $N$ Galois representation attached to an elliptic curve $E/\bbQ$, then the determinant character is surjective - i.e., $\det(H) = (\bbZ/N\bbZ)^\times$. We assume from here on that $\det(H) = (\bbZ/N\bbZ)^\times$ (in particular this is the case for Cartan subgroups of $\GL_2(\bbZ/p^l\bbZ)$) so that $X(H)$ consists of one geometrically irreducible component. 

			Suppose that $N$ factors as a product $p_1^{l_1}...p_r^{l_r}$, where $p_1, ..., p_r$ are distinct primes. Let 
			$$H = H_1 \times ... \times H_r \subset \GL_2(\bbZ/p_1^{l_1}\bbZ) \times ... \times \GL_2(\bbZ/p_r^{l_r}\bbZ) \cong \GL_2(\bbZ/N\bbZ)$$ 
			be a subgroup (such that $-I \in H$ and $\det(H) = (\bbZ/N\bbZ)^\times$). In this case we denote the curve $X(H)$ by $X(H_1, ..., H_r)$. When the groups $H_i$ are equal to (a conjugate of) $G(\s  p^l)$, $G(\s  {p^l}^+)$, $G(\ns  p^l)$, or $G(\ns  {p^l}^+)$ we will replace $H_i$ with the label $\s  p^l$, $\s  {p^l}^+$, $\ns  p^l$, and $\ns  {p^l}^+$ respectively.
		
	\section{Quadratic Twist-Type \texorpdfstring{$p^l$}{p{\textasciicircum}l}-Congruences of Elliptic Curves} \label{sec:technical}
		
		The idea for the following lemma is due to Cremona--Freitas \cite[Lemma 2.2 and Theorem 2.4]{CF_GMFTSTOCBEC} and Halberstadt \cite[Section 2.1]{H_SLCMX11}. We provide a generalisation to composite $N$.

		Let $K$ be a field of characteristic zero. For an elliptic curve $E/K$ we let $E^d/K$ denote the quadratic twist of $E$ by $d \in K^\times$.
		
		\begin{lemma} \thlabel{cong_q_twist}
			Let $E/K$ be an elliptic curve and $N > 2$ be an integer. Let $H^+ = \bar{\rho}_{E, N}(G_K)$ be the image of the mod $N$ Galois representation $\bar{\rho}_{E,N}$. Then $E$ admits a quadratic twist-type $N$-congruence if and only if there exists an index $2$ subgroup $H \subset H^+$ and a matrix $g \in \GL_2(\bbZ/N\bbZ)$ such that 
			\begin{align} \label{cong_cond}
				g^{-1}  h g = h \quad \textnormal{for all $h \in H$} \qquad \textnormal{and} \qquad
				g^{-1} hg = -h \quad \textnormal{for all $h \in H^+ \setminus H$}.
			\end{align}
			More explicitly, let $d \in K^\times$ be such that $\bar{\rho}_{E, N}^{-1}(H) = G_{K(\sqrt{d})}$. Then $E$ is $N$-congruent to $E^d$ and the power of the congruence is $\det(g)$.
			
			\begin{proof}
				For any $d \in K^\times$ let $t \colon E \to E^d$ be an isomorphism defined over $K(\sqrt{d})$. Then for any $P \in E$ and any $\sigma \in G_K$ we have  
				\begin{equation*}
					\sigma(t(P)) = \varepsilon_d(\sigma) t(\sigma(P))
				\end{equation*}
				where $\varepsilon_d \colon G_K \to \{\pm 1\}$ is the quadratic character attached to $K(\sqrt{d})/K$. Fix a basis for $E[N]$ and use $t$ to choose a compatible basis for $E^d[N]$. Then for any $\sigma \in G_K$ we have that 
				\begin{equation} \label{rep_q_twist}
					\bar{\rho}_{E^d, N}(\sigma) = \varepsilon_d(\sigma)\bar{\rho}_{E, N}(\sigma).
				\end{equation}
				
				Suppose that $\phi \colon E^d[N] \to E[N]$ is an $N$-congruence between $E$ and its quadratic twist $E^d$. Then let $\phi$ be given by $g \in \GL_2(\bbZ/N\bbZ)$ with respect to the chosen bases for $E[N]$ and $E^d[N]$. Note that $g^{-1} \bar{\rho}_{E, N}(\sigma) g = \bar{\rho}_{E^d, N}(\sigma)$ for all $\sigma \in G_K$. Let $H = \bar{\rho}_{E, N}(G_{K(\sqrt{d})})$. By substituting into (\ref{rep_q_twist}) we see that either (\ref{cong_cond}) holds, or $H = H^+$.
				
				Suppose that $H = H^+$. Since $d$ is not a square, there exists $\sigma, \tau \in G_K$ such that $\bar{\rho}_{E, N}(\sigma) = \bar{\rho}_{E, N}(\tau)$ but $\varepsilon_d(\sigma) = -1$ and $\varepsilon_d(\tau) = 1$. Then by (\ref{rep_q_twist}) we have 
				\begin{equation*}
					g^{-1} \bar{\rho}_{E, N}(\sigma) g = -\bar{\rho}_{E, N}(\sigma) = -\bar{\rho}_{E, N}(\tau) = - g^{-1} \bar{\rho}_{E, N}(\tau) g.
				\end{equation*}
				So $2 \bar{\rho}_{E, N}(\sigma) = 0$, contradicting the assumption that $N > 2$. Hence $H$ has index $2$ in $H^+$.
				
				Conversely, take $d \in K^\times$ so that $G_{K(\sqrt{d})} = \bar{\rho}_{E, N}^{-1}(H)$. Then by substituting (\ref{cong_cond}) into (\ref{rep_q_twist}) we see that
				\begin{align*} 
					g^{-1} \bar{\rho}_{E, N}(\sigma) g = \varepsilon_d(\sigma)\bar{\rho}_{E, N}(\sigma) = \bar{\rho}_{E^d, N}(\sigma) \quad \textnormal{for all $\sigma \in G_K$}.
				\end{align*}
				Thus $g$ gives an $N$-congruence between $E$ and $E^d$.
				
				Finally note that $t|_{E[N]}$ preserves the Weil pairing. Therefore
				\begin{equation*}
					e_{E^d, N}(\phi^{-1}(P), \phi^{-1}(Q)) = e_{E^d, N}(t(P), t(Q))^{\det(g)} = e_{E, N}(P, Q)^{\det(g)}.
				\end{equation*}
				Hence $\phi$ has power $\det(g)$, as required.
			\end{proof}
		\end{lemma}

		\begin{remark}
			Note that the assumption that $N > 2$ in \thref{cong_q_twist} is necessary. An elliptic curve $E/K$ is $2$-congruent to every quadratic twist, so there may exist $d \in K^\times$ such that $E$ is $2$-congruent to $E^d$ but $\bar{\rho}_{E, 2}(G_K) =  \bar{\rho}_{E, 2}(G_{K(\sqrt{d})})$.
		\end{remark}

		\begin{remark} \thlabel{rem:traces}
			Taking traces in (\ref{cong_cond}) shows that for any $h \in H^+ \setminus H$ we have $2\Tr(h) = 0$ and $2\Tr(g) = 0$. Hence $H$ contains the subgroup of $H^+$ generated by $\{h \in H^+ : 2\Tr(h) \neq 0 \}$. 

			When $N = 2^k$ is a power of $2$ this is computationally useful. Since $| \bar{\rho}_{E, N}(G_K) |$ may be divisible by a high power of $2$ there may be many index $2$ subgroups to check. However noting that $H$ contains the subgroup generated by $\{h \in H^+ : 2\Tr(h) \neq 0 \}$ significantly reduces the number of subgroups to consider.

			For a given subgroup $H^+ \subset \GL_2(\bbZ/N\bbZ)$ our \texttt{Magma} function \texttt{GL2GivesQTTCongruence} uses this observation to determine all such $H$, $g$ satisfying (\ref{cong_cond}) (see \cite{ME_ELECTRONIC}). 
		\end{remark}

		\subsection{Quadratic Twist-Type \texorpdfstring{$p^l$}{p{\textasciicircum}l}-Congruences}
	
		\thref{cong_q_twist} allows us to extend \cite[Theorem 2.4]{CF_GMFTSTOCBEC} to powers of odd primes (while also ensuring we avoid their condition \textbf{(S)}). For an odd prime $p$ let $\xi$ be a quadratic nonresidue modulo $p$.
		
		\begin{prop} \thlabel{cartan_cong_pl}
			
			Let $E/K$ be an elliptic curve, let $p$ be an odd prime, and let $l \geq 1$ be an integer. Then $E$ admits a quadratic twist-type $(p^l, -1)$-congruence if and only if there exists a choice of basis for $E[p^l]$ such that the image of the mod $p^l$ Galois representation, $\bar{\rho}_{E, p^l}(G_K)$, is contained in $G(\s {p^l}^+)$, but not the subgroup $G(\s p^l)$. More explicitly, let $d \in K^\times$ be such that $\bar{\rho}_{E, p^l}(G_{K(\sqrt{d})}) = \bar{\rho}_{E, p^l}(G_K) \cap G(\s {p^l})$. Then $E$ is $(p^l, -1)$-congruent to $E^d$.
			
			Similarly $E$ admits a quadratic twist-type $(p^l, -\xi)$-congruence if and only if there exists a choice of basis for $E[p^l]$ such that the image of the mod $p^l$ Galois representation, $\bar{\rho}_{E, p^l}(G_K)$, is contained in $G(\ns {p^l}^+)$, but not the subgroup $G(\ns p^l)$. More explicitly, let $d \in K^\times$ be such that $\bar{\rho}_{E, p^l}(G_{K(\sqrt{d})}) = \bar{\rho}_{E, p^l}(G_K) \cap G(\ns {p^l})$. Then $E$ is $(p^l, -\xi)$-congruent to $E^d$.
							
			\begin{proof}
				First suppose that the image of the mod $p^l$  Galois representation, $\bar{\rho}_{E, p^l}(G_K)$, is contained in $G(\s {p^l}^+)$, and $\bar{\rho}_{E, p^l}(G_{K(\sqrt{d})}) = \bar{\rho}_{E, p^l}(G_K) \cap G(\s {p^l})$. 
				
				By \thref{cong_q_twist} the elliptic curve $E$ is $(p^l, -1)$-congruent to $E^d$ if there exists a matrix $g \in \GL_2(\bbZ/p^l\bbZ)$ which commutes with the elements of $G(\s  {p^l})$, and such that $g^{-1} \big(\begin{smallmatrix} 0 & 1\\ 1 & 0 \end{smallmatrix}\big) g = -\big(\begin{smallmatrix} 0 & 1\\ 1 & 0 \end{smallmatrix}\big)$. 
				
				We define the matrix $g_{\s} = \big(\begin{smallmatrix} 1 & 0\\ 0 & -1 \end{smallmatrix}\big) \in G(\s  p^l)$ and we claim $g_{\s}$ satisfies these conditions. Now $g_{\s}$ commutes with the elements of $G(\s  {p^l})$ since $G(\s  {p^l})$ is abelian and $g_{\s} \in G(\s p^l)$. We may directly verify that $g_{\s}^{-1} \big(\begin{smallmatrix} 0 & 1\\ 1 & 0 \end{smallmatrix}\big) g_{\s} = -\big(\begin{smallmatrix} 0 & 1\\ 1 & 0 \end{smallmatrix}\big)$.
				
				Similarly in the nonsplit case suppose that the image of the mod $p^l$  Galois representation, $\bar{\rho}_{E, p^l}(G_K)$, is contained in $G(\ns {p^l}^+)$, and $\bar{\rho}_{E, p^l}(G_{K(\sqrt{d})}) = \bar{\rho}_{E, p^l}(G_K) \cap G(\ns {p^l})$. 
				
				We define the matrix $g_{\ns} = \big(\begin{smallmatrix} 0 & \xi \\ 1 & 0 \end{smallmatrix}\big)$. Now $g_{\ns}$ commutes with the elements of $G(\ns  {p^l})$ since $G(\ns  {p^l})$ is abelian and $g_{\ns} \in G(\ns p^l)$. A direct check verifies $g_{\ns}^{-1} \big(\begin{smallmatrix} 1 & 0\\ 0 & -1 \end{smallmatrix}\big) g_{\ns} = -\big(\begin{smallmatrix} 1 & 0\\ 0 & -1 \end{smallmatrix}\big)$. It follows from \thref{cong_q_twist} that $E$ is $(p^l, -\xi)$-congruent to $E^d$.

				For the converse suppose that $E$ is $p^l$-congruent (over $K$) to a nontrivial quadratic twist, $E^d$. Let $H = \bar{\rho}_{E, p^l}(G_{K(\sqrt{d})})$, $H^+ = \bar{\rho}_{E, p^l}(G_{K})$ and $g$ be as given by \thref{cong_q_twist}.

				As noted in \thref{rem:traces} we have $\Tr(g) = 0$. In particular by the Cayley--Hamilton Theorem the image of $g$ in $\PGL_2(\bbZ/p^l \bbZ)$ has order $2$. The following lemma places a strong restriction on the possible $g$.

				\begin{lemma} \thlabel{lemma:order2pgl}
					There are exactly two conjugacy classes of elements of order $2$ in $\PGL_2(\bbZ/p^l\bbZ)$. More precisely if the image of an element $g \in \GL_2(\bbZ/p^l\bbZ)$ in $\PGL_2(\bbZ/p^l\bbZ)$ has order $2$, then $g$ is conjugate (in $\PGL_2(\bbZ/p^l\bbZ)$) to $g_{\s} = \big(\begin{smallmatrix} 1 & 0\\ 0 & -1 \end{smallmatrix}\big)$ if $-\det(g)$ is a quadratic residue modulo $p$ or $g_{\ns} = \big(\begin{smallmatrix} 0 & \xi \\ 1 & 0 \end{smallmatrix}\big)$ if $-\det(g)$ is a quadratic nonresidue modulo $p$.

					\begin{proof}
						By the Cayley--Hamilton Theorem an element $g \in \PGL_2(\bbZ/p^l\bbZ)$ has order $2$ if and only if $\Tr(g) = 0$. We may then write $g = \big(\begin{smallmatrix} a & b\\ c & -a \end{smallmatrix}\big)$.

						Now if $-\det(g)$ is a quadratic residue modulo $p$ let $\beta \in \bbZ/p^l\bbZ$ be such that $\beta^2 = -\det(g)$. Consider $\gamma = \big(\begin{smallmatrix} \beta + a & -\beta + a\\ c & c \end{smallmatrix}\big)$ and note that $\gamma^{-1} g \gamma = \beta g_{\s}$.
						
						Similarly if $-\det(g)$ is a quadratic nonresidue modulo $p$ let $\beta \in \bbZ/p^l\bbZ$ be such that $\beta^2 = -\det(g)/\xi$. Then consider $\gamma = \big(\begin{smallmatrix} a + \beta & \xi\beta + a\\ c & c \end{smallmatrix}\big)$ and note that $\gamma^{-1} g \gamma = \beta g_{\ns}$.
					\end{proof}
				\end{lemma}
				
				By \thref{lemma:order2pgl} the matrix $g$ is conjugate to either $g_{\s}$ or $g_{\ns}$ in $\PGL_2(\bbZ/p^l\bbZ)$ depending on whether the congruence given by $g$ has power $-1$ or $-\xi$ respectively.
				
				In the first case (i.e., if $E$ is $(p^l, -1)$-congruent to $E^d$) we are free to choose a basis for $E[p^l]$ so that $g = u g_{\s}$ for some scalar $u \in (\bbZ/p^l\bbZ)^\times$. But by \thref{cong_q_twist} the subgroup $H$ is contained in the centraliser of $g$ in $\GL_2(\bbZ/p^l\bbZ)$. A direct calculation shows that $H$ is contained in $G(\s p^l)$. By considering matrices $h \in \GL_2(\bbZ/p^l\bbZ)$ such that $g h = -h g$ we see that $H^+$ is contained in $G(\s {p^l}^+)$ and $H = H^+ \cap G(\s p^l)$. 
				
				In the second case (i.e., if $E$ is $(p^l, -\xi)$-congruent to $E^d$) we are free to choose a basis for $E[p^l]$ so that $g = u g_{\ns}$ for some scalar $u \in (\bbZ/p^l\bbZ)^\times$. But by \thref{cong_q_twist} the subgroup $H$ is contained in the centraliser of $g$ in $\GL_2(\bbZ/p^l\bbZ)$, so $H$ is contained in $G(\ns p^l)$. By considering matrices $h \in \GL_2(\bbZ/p^l\bbZ)$ such that $g h = -h g$ we see that $H^+$ is contained in $G(\ns {p^l}^+)$ and $H = H^+ \cap G(\s p^l)$.
				
			\end{proof}
		\end{prop}

	\subsection{Quadratic Twist-Type \texorpdfstring{$2^k$}{2{\textasciicircum}k}-Congruences} \label{sec:qtt2}
		
		It is clear that an elliptic curve $E/K$ is $2$-congruent to any quadratic twist $E^d/K$ for $d \in K^\times$. Furthermore, it was shown in \cite[Corollary 7.4]{BD_TAOG2CW44SJ} that any elliptic curve is $(4,3)$-congruent to its quadratic twist by its discriminant. Taking $k = 1$ in the following proposition removes the hypothesis on the image of $G_K$ (since $G(\ns  2^+) = \GL_2(\bbZ/2\bbZ)$) and we recover \cite[Corollary 7.4]{BD_TAOG2CW44SJ}. In fact, we also recover their proof.
		
		\begin{prop} \thlabel{2r_cong_disc}
			Let $E/K$ be an elliptic curve and let $k \geq 1$ be an integer. Suppose that there exists a choice of basis for $E[2^k]$ such that the image of the mod $2^k$ Galois representation, $\bar{\rho}_{E, 2^{k}}(G_K)$, is contained in the subgroup $G(\ns {2^k}^+) \subset \GL_2(\bbZ/2^k\bbZ)$. Then $E$ is $(2^{k+1}, 3)$-congruent to its quadratic twist by its discriminant $\Delta \colonequals \Delta(E)$.

			Similarly, if there exists a choice of basis for $E[2^k]$ such that the image of the mod $2^k$ Galois representation, $\bar{\rho}_{E, 2^{k}}(G_K)$, is contained in the subgroup\footnote{Recall from Section \ref{sec:basic_defs} that $G(\s {2^k}^+)$ is defined to be the subgroup generated by $G(\s 2^k)$ and $\big(\begin{smallmatrix} 0 & 1\\ 1 & 0 \end{smallmatrix}\big)$ \emph{not} as the normaliser of $G(\s {2^k})$ in $\GL_2(\bbZ/2^k\bbZ)$.} $G(\s {2^k}^+) \subset \GL_2(\bbZ/2^k\bbZ)$ then $E$ is $(2^{k+1}, -1)$-congruent to its quadratic twist by its discriminant.
			
			\begin{proof}
				Let $\pi_{2^k} \colon \GL_2(\bbZ/2^{k+1}\bbZ) \to \GL_2(\bbZ/2^k\bbZ)$ be the natural projection. 

				First note that if $E$ is given by a Weierstrass equation $y^2 = f(x)$ then $\Delta$ is equal to the discriminant of $f(x)$ up to a square factor. Now $G(\ns {2})$ is the unique index $2$ subgroup of $\GL_2(\bbZ/2\bbZ)$, so $\bar{\rho}_{E, 2}(G_{K(\sqrt{\Delta})}) = G(\ns {2}) \cap \bar{\rho}_{E, 2}(G_K)$. For any $k \geq 1$ the image of $G(\ns {2^k})$ under the projection $\GL_2(\bbZ/2^{k}\bbZ) \to \GL_2(\bbZ/2\bbZ)$ is $G(\ns {2})$. Therefore if $\bar{\rho}_{E, 2^{k}}(G_K)$ is contained in $G(\ns {2^k})$ the discriminant, $\Delta$, is a square in $K$.
	
				Combining this observation with \thref{cong_q_twist} we see that in the nonsplit case it suffices to find a matrix $g \in \GL_2(\bbZ/2^{k+1}\bbZ)$ with determinant $3$ which commutes with the elements of $\pi_{2^k}^{-1}G(\ns {2^k})$, and such that $g^{-1}\big(\begin{smallmatrix} 1 & -1\\ 0 & -1 \end{smallmatrix}\big)g = -\big(\begin{smallmatrix} 1 & -1\\ 0 & -1 \end{smallmatrix}\big)$. Note that if $\bar{\rho}_{E, 2^{k}}(G_K)$ is contained in $G(\ns {2^k})$ then a similar argument to \thref{cong_q_twist} shows that such a $g$ induces a $G_K$-equivariant automorphism of $E[2^{k+1}]$ with power $3$. We claim that $g = \big(\begin{smallmatrix} 1 & -2\\ 2 & -1 \end{smallmatrix}\big)$ satisfies this condition. 
				
				A direct check verifies that $g^{-1}\big(\begin{smallmatrix} 1 & -1\\ 0 & -1 \end{smallmatrix}\big)g = -\big(\begin{smallmatrix} 1 & -1\\ 0 & -1 \end{smallmatrix}\big)$. Now $\pi_{2^k}^{-1}G(\ns {2^k})$ contains $G(\ns {2^{k+1}})$ as an abelian subgroup of index $4$. Moreover, $\pi_{2^k}^{-1}G(\ns {2^k})$ is generated by $G(\ns {2^{k+1}})$ and the matrices of the form $I + 2^k M$. But since $g$ is contained in $G(\ns {2^{k+1}})$ it commutes with the elements of $G(\ns {2^{k+1}})$. Note that $g^{-1}(I + 2^k M)g = I + 2^k g^{-1} M g$ depends only on $g \equiv I$ modulo $2$. Thus $g$ commutes with the elements of $\pi_{2^k}^{-1}G(\ns {2^k})$.

				Similarly in the split case by \thref{cong_q_twist} it suffices to find a matrix $g \in \GL_2(\bbZ/2^{k+1}\bbZ)$ with determinant $-1$ which commutes with the elements of $\pi_{2^k}^{-1}G(\s {2^k})$, and such that $g^{-1}\big(\begin{smallmatrix} 0 & 1\\ 1 & 0 \end{smallmatrix}\big)g = -\big(\begin{smallmatrix} 0 & 1\\ 1 & 0 \end{smallmatrix}\big)$. We claim that $g = \big(\begin{smallmatrix} 1 & 0\\ 0 & -1 \end{smallmatrix}\big)$ satisfies this condition. 
				
				A direct check verifies that $g^{-1}\big(\begin{smallmatrix} 0 & 1\\ 1 & 0 \end{smallmatrix}\big)g = -\big(\begin{smallmatrix} 0 & 1\\ 1 & 0 \end{smallmatrix}\big)$. Since $g \in G(\s {2^{k+1}})$ it commutes with the elements of $G(\s {2^{k+1}})$, and as in the nonsplit case $g^{-1}(I + 2^k M)g = I + 2^k M$ -- so $g$ commutes with elements of $\pi_{2^k}^{-1}G(\s {2^k})$. 
			\end{proof}
		\end{prop}
		
		Now suppose $K = \bbQ$. We wish to give a classification of quadratic twist-type $(2^k, r)$-congruences of elliptic curves over $\bbQ$ in terms of the image of their mod $2^k$ Galois representations, mimicking that given by \thref{cartan_cong_pl} in the case of odd prime powers.
		
		By analysing rational points on modular curves of level $2^k$, Rouse and Zureick-Brown \cite{RZB_ECOQA2AIOG} have classified all possible images of the $2$-adic Galois representations attached to elliptic curves $E/\bbQ$ without CM. Specifically Rouse and Zureick-Brown determine a list of $1208$ open subgroups $G_i \subset \GL_2(\bbZ_2)$ for which there exists an elliptic curve $E/\bbQ$ without CM such that the image of the $2$-adic Galois representation, $\rho_{E, 2}(G_{\bbQ})$, is conjugate to $G_i$. 
		
		It will be more convenient for us to refer to subgroups of $\GL_2(\bbZ/2^k\bbZ)$ than to their inverse images in $\GL_2(\bbZ_2)$ and in Table \ref{table:2adiccongs} we introduce $11$ subgroups of $\GL_2(\bbZ/2^k\bbZ)$ which give rise to \emph{all} quadratic twist-type $(2^k, r)$-congruences of non-CM elliptic curves over $\bbQ$. The labelling we use for these subgroups corresponds to that of Rouse and Zureick-Brown in the following way: a subgroup $G_i \subset \GL_2(\bbZ_2)$ labelled by Rouse--Zureick-Brown is equal to the inverse image under the projection $\GL_2(\bbZ_2) \to \GL_2(\bbZ/2^k\bbZ)$ of the subgroup we label $G_i$.

		It will often be convenient for us to refer to such an index $2$ subgroup $H \subset H^+ = \bar{\rho}_{E, 2^k}(G_K)$ in terms of $H^+$ and $E$. We therefore define the following notation. Let $d \in K^\times$ be such that $K(\sqrt{d})/K$ is the quadratic extension cut out by $\bar{\rho}_{E, 2^k}^{-1}(H)$. Then we denote $H$ by $H^+(\sqrt{d})$ (note that $d$ may depend on $E$).
		
		Since \thref{cong_q_twist} is stated purely in terms of the mod $N$ Galois representation of an elliptic curve, the classification of Rouse and Zureick-Brown allows us to give a complete description of \emph{non-CM} elliptic curves defined over $\bbQ$ admitting a $(2^k, r)$-congruence of quadratic twist-type. 
		
		In particular we have:

		\begin{prop} \thlabel{prop:2adiccongs}
			Let $E/\bbQ$ be an elliptic curve without CM and let $k \geq 2$ be an integer. Then $E$ is $(2^k,r)$-congruent to a nontrivial quadratic twist $E^d$ if and only if there exists a choice of basis for $E[2^k]$ such that the image of the mod $2^k$ Galois representation, $\bar{\rho}_{E, 2^{k}}(G_{\bbQ})$, is contained in one of the subgroups $H^+ \subset \GL_2(\bbZ/2^k\bbZ)$ listed in Table \ref{table:2adiccongs}, but not the corresponding subgroup $H$, and $\bar{\rho}_{E, 2^k}(G_{\bbQ(\sqrt{d})}) = \bar{\rho}_{E, 2^k}(G_{\bbQ}) \cap H$. 
		\end{prop}

\begin{table}[ht]
			\centering
			\begin{adjustbox}{width=\textwidth}

			\begin{tabular}{c|c|c|c|c|c|c}
\multirow{2}{*}{$(2^k,r)$ } & \multicolumn{3}{c|}{The subgroup $H^+$ of $\GL_2(\bbZ/2^k\bbZ)$} & \multicolumn{2}{c|}{The index $2$ subgroup $H \subset H^+$} & \multirow{2}{*}{Quadratic Twist}  \\

						& Generators 	& Index 	& Label 						& Generators						& Label		&				\\ 
\hline
\multirow{2}{*}{$(4,1)$}

						& {$\left\{ \big(\begin{smallmatrix} 1 & 1\\ 0 & 1 \end{smallmatrix}\big) , \big(\begin{smallmatrix} 3 & 0\\ 2 & 1 \end{smallmatrix}\big) \right\}$}
															& $6$	&{$G_{9}$}	& $\left\{ \big(\begin{smallmatrix} 3 & 0\\ 0 & 3 \end{smallmatrix}\big) , \big(\begin{smallmatrix} 1 & 1\\ 0 								& 1 \end{smallmatrix}\big) \right\}$
																						 									& $G_{9}(\sqrt{-1})$	& $-1$			\\

						& {$\left\{ \big(\begin{smallmatrix} 1 & 0\\ 2 & 3 \end{smallmatrix}\big) , \big(\begin{smallmatrix} 1 & 2\\ 0 & 1 \end{smallmatrix}\big), \big(\begin{smallmatrix} 3 & 1\\ 2 & 1 \end{smallmatrix}\big) \right\}$}
															& $6$	&{$G_{11}$}	& $\left\{ \big(\begin{smallmatrix} 1 & 2\\ 0 & 1 \end{smallmatrix}\big) , \big(\begin{smallmatrix} 1 & 3\\ 2 								& 3 \end{smallmatrix}\big) \right\}$									
																															& $G_{11}(\sqrt{-1})$	& $-1$			\\
\hdashline
\multirow{3}{*}{$(4,3)$}
						& $\GL_2(\bbZ/4\bbZ)$				& $1$	& $G_{1}$ 					& $\pi_2^{-1} G(\ns 2)$				& $-$	& $\Delta(E)$	\\

						& {$\left\{ \big(\begin{smallmatrix} 1 & 1\\ 0 & 1 \end{smallmatrix}\big) , \big(\begin{smallmatrix} 3 & 0\\ 2 & 1 \end{smallmatrix}\big) \right\}$}
															& $6$	&{$G_{9}$}					& $\left\{ \big(\begin{smallmatrix} 1 & 2\\ 0 & 1 \end{smallmatrix}\big) , \big(\begin{smallmatrix} 1 & 1\\ 2 & 1 \end{smallmatrix}\big) \right\}$
																																	& $G_{9}(\sqrt{-\Delta})$	& $-\Delta(E)$	\\

						& {$\left\{ \big(\begin{smallmatrix} 1 & 0\\ 2 & 3 \end{smallmatrix}\big) , \big(\begin{smallmatrix} 1 & 2\\ 0 & 1 \end{smallmatrix}\big), \big(\begin{smallmatrix} 3 & 1\\ 2 & 1 \end{smallmatrix}\big) \right\}$}
															& $6$ 	&{$G_{11}$}	& $\left\{ \big(\begin{smallmatrix} 1 & 2\\ 0 & 1 \end{smallmatrix}\big) , \big(\begin{smallmatrix} 3 & 1\\ 0 																	& 1 \end{smallmatrix}\big), \big(\begin{smallmatrix} 3 & 0\\ 0 & 3 \end{smallmatrix}\big) \right\}$	
																															& $G_{11}(\sqrt{-\Delta})$ 	& $-\Delta(E)$ 	\\
\hline
$(8,1)$					& $\left\{ \big(\begin{smallmatrix} 1 & 3\\ 0 & 7 \end{smallmatrix}\big) , \big(\begin{smallmatrix} 1 & 3\\ 2 & 7 \end{smallmatrix}\big), \big(\begin{smallmatrix} 1 & 1\\ 6 & 3 \end{smallmatrix}\big) \right\}$
						& $24$ 	& $G_{97}$					& $\left\{ \big(\begin{smallmatrix} 3 & 4\\ 0 & 3 \end{smallmatrix}\big) , \big(\begin{smallmatrix} 5 & 6\\ 4 								& 1 \end{smallmatrix}\big), \big(\begin{smallmatrix} 5 & 1\\ 6 & 7 \end{smallmatrix}\big) \right\}$
																						& $G_{97}(\sqrt{-1})$	& $-1$			\\ 
\hdashline
\multirow{2}{*}{$(8,3)$}   	
						& $\pi_4^{-1} G(\ns 4^+)$			& $4$	&$G_{7}$					& $\pi_4^{-1} G(\ns 4)$				& $-$ 	& $\Delta(E)$	\\

						& $\left\{\big(\begin{smallmatrix} 1 & 3\\ 6 & 7 \end{smallmatrix}\big) , \big(\begin{smallmatrix} 1 & 1\\ 4 & 3 \end{smallmatrix}\big), \big(\begin{smallmatrix} 1 & 2\\ 0 & 5 \end{smallmatrix}\big) \right\}$ 
						& $24$	& $G_{90}$					& $\left\{\big(\begin{smallmatrix} 3 & 4\\ 0 & 3 \end{smallmatrix}\big) , \big(\begin{smallmatrix} 3 & 3\\ 4 								& 1 \end{smallmatrix}\big), \big(\begin{smallmatrix} 5 & 6\\ 0 & 1 \end{smallmatrix}\big) \right\}$
																						& $G_{90}(\sqrt{-\Delta})$	& $-\Delta(E)$ 	\\
\hdashline
$(8,5)$					& $\left\{\big(\begin{smallmatrix} 1 & 2\\ 2 & 7 \end{smallmatrix}\big) , \big(\begin{smallmatrix} 1 & 3\\ 2 & 7 \end{smallmatrix}\big), \big(\begin{smallmatrix} 1 & 1\\ 6 & 3 \end{smallmatrix}\big) \right\}$                                 
						& $24$	& $G_{63}$					& $\left\{\big(\begin{smallmatrix} 3 & 4\\ 0 & 3 \end{smallmatrix}\big) , \big(\begin{smallmatrix} 5 & 6\\ 4 								& 1 \end{smallmatrix}\big), \big(\begin{smallmatrix} 5 & 1\\ 6 & 7 \end{smallmatrix}\big) \right\}$
																						& $G_{63}(\sqrt{-1})$ 	& $-1$ 			\\

\hdashline
\multirow{2}{*}{$(8,7)$}	
						& $\pi_4^{-1} G(\s 4^+)$     	    & $12$ &$G_{23}$					& $\pi_4^{-1} G(\s 4)$				& $-$ 	& $\Delta(E)$ 	\\

						& $\left\{ \big(\begin{smallmatrix} 1 & 3\\ 0 & 7 \end{smallmatrix}\big) , \big(\begin{smallmatrix} 1 & 1\\ 2 & 7 \end{smallmatrix}\big) , \big(\begin{smallmatrix} 1 & 2\\ 0 & 5 \end{smallmatrix}\big) \right\}$ 
															& $24$	&$G_{91}$ 					& $\left\{ \big(\begin{smallmatrix} 3 & 4\\ 0 & 3 \end{smallmatrix}\big) , \big(\begin{smallmatrix} 1 & 7\\ 0 								& 7 \end{smallmatrix}\big), \big(\begin{smallmatrix} 7 & 1\\ 0 & 1 \end{smallmatrix}\big), \big										(\begin{smallmatrix} 5 & 6\\ 0 & 1 \end{smallmatrix}\big) \right\} $
																															& $G_{91}(\sqrt{-\Delta})$	& $-\Delta(E)$ 	\\

\hline
$(16,3)$				& $\pi_8^{-1} G(\ns 8^+)$			& $16$	&$G_{55}$					& $\pi_8^{-1} G(\ns 8)$ 			& $-$	 	& $\Delta(E)$ 	\\ 
\hline
$(32,3)$				& $\pi_{16}^{-1} G(\ns 16^+)$		& $64$	& $G_{441}$					& $\pi_{16}^{-1} G(\ns 16)$			& $-$	 	& $\Delta(E)$			
				
			\end{tabular}
			
			\end{adjustbox}

			\caption{Subgroups of $\GL_2(\bbZ/2^k\bbZ)$ which give rise to quadratic twist-type $(2^k, r)$-congruences over $\bbQ$. Here $\pi_{2^k}$ denotes the natural reduction map $\GL_2(\bbZ/2^{k+1}\bbZ) \to \GL_2(\bbZ/2^{k}\bbZ)$. As noted in \thref{remark:conjugatetotranspose}, the subgroups $H$ and $H^+$ above are conjugate to their transposes (this is checked in the electronic data \cite{ME_ELECTRONIC}). The modular curves $X(H^+)$ all have genus $0$ except $X(\ns16^+)$ which has genus $2$. } 
			\label{table:2adiccongs}
\end{table}

			\begin{proof}
				This is a brute force calculation in \texttt{Magma} using the classification of Rouse and Zureick-Brown \cite{RZB_ECOQA2AIOG} and \thref{cong_q_twist}. The observations in \thref{rem:traces} significantly speed up the search.
			\end{proof}

		\begin{example}[A $32$-Congruence] \thlabel{ex:32}
			Heegner \cite{H_DAUM}, Baran \cite{B_NONSCSMCATCNP} and Rouse--Zureick-Brown \cite{RZB_ECOQA2AIOG} determined the rational points on the genus $2$ curve $X(\ns  16^+)$ (see \cite[Remark 1.5]{RZB_ECOQA2AIOG} for a discussion). There are \emph{exactly} two non-CM $j$-invariants for which there are rational points in the fibre of $X(\ns  16^+) \to X(1)$ above $j$, namely:
			\begin{equation*}
				\frac{-2^{18} \cdot 3 \cdot 5^3 \cdot 13^3 \cdot 41^3 \cdot 107^3}{17^{16}} \qquad \textnormal{and} \qquad \frac{-2^{21} \cdot 3^3 \cdot 5^3 \cdot 7 \cdot 13^3 \cdot 23^3 \cdot 41^3 \cdot 179^3 \cdot 409^3}{79^{16}}.
			\end{equation*}
			
			The elliptic curves $E_{(32,3)a}/\bbQ$ and $E_{(32,3)b}/\bbQ$ from the statement of \thref{thm:examples} are the minimal twists with the $j$-invariants above. By \thref{2r_cong_disc} the elliptic curves $E_{(32,3)a}$ and $E_{(32,3)b}$ admit $(32, 3)$-congruences with their quadratic twists by their discriminants. In fact, these are the only $j$-invariants of elliptic curves admitting nontrivial quadratic twist-type $32$-congruences over $\bbQ$ (see \thref{thm:infmanypl}).
		\end{example}
	
	\section{Trivial Quadratic Twist-Type Congruences and CM} \label{sec:trivialCM}

		Let $p$ be an odd prime and let $l \geq 1$ be an integer. Recall that if $E/\bbQ$ has CM by an imaginary quadratic order $\mathcal{O}$ of discriminant $-D$ and $p$ is unramified in $\mathcal{O}$ then $\bar{\rho}_{E, p^l}(G_{\bbQ})$ is contained in the normaliser of a Cartan subgroup of $\GL_2(\bbZ/p^l\bbZ)$ and the inverse image of the corresponding Cartan subgroup is $G_{\bbQ(\sqrt{-D})}$. Moreover the Cartan subgroup is nonsplit if and only if $p$ is inert in $\mathbb{Q}(\sqrt{-D})$ (see \cite[A.5]{S_LOTMWT} or \cite[Proposition 4.1]{B_NONSCSMCATCNP}). It then follows from \thref{cartan_cong_pl} that CM elliptic curves admit quadratic twist-type $p^l$-congruences for all odd prime numbers $p$ unramified in the CM field. In fact, quadratic twist-type congruences are trivial if and only if they arise in this way:

		\begin{lemma} \thlabel{lemma:CM_congs}
			An elliptic curve $E/\bbQ$ admits an isogeny with a nontrivial quadratic twist if and only if $E$ has CM. More precisely if $E$ has CM by an imaginary quadratic order, $\mathcal{O}$, of discriminant $-D$ then $E$ is $D$-isogenous to its quadratic twist by $-D$. Moreover if $4$ divides $D$ then we may take the isogeny to have degree $D/4$.

			\begin{proof}
				Suppose that an elliptic curve $E/\bbQ$ admits an isogeny with a nontrivial quadratic twist $E^d/\bbQ$ (where $d \in \mathbb{Z}$ is possibly positive). Composing the isogeny with the $\bbQ(\sqrt{d})$-isomorphism $E^d \to E$ gives rise to an endomorphism of $E$ defined over $\bbQ(\sqrt{d})$ but not over $\bbQ$ (hence is not multiplication by $m$). Thus $E$ has CM by an order in $\bbQ(\sqrt{d})$ (by \cite[Theorem VI.6.1(b)]{S_TAOEC}).

				Conversely, if $4$ divides $D$ (i.e., if $-D \in \{-4, -8, -12, -16, -28\}$) let $D' = D/4$ and let $D' = D$ otherwise. Note that $\sqrt{-D'} \in \mathcal{O}$. Composing the $\bbQ(\sqrt{-D'})$-isomorphism $E \to E^{-D'}$ with the endomorphism $\sqrt{-D'}$ gives an isogeny $E \to E^{-D'}$ defined over $\bbQ$. The degree of the isogeny is $D'$, since the endomorphism $\sqrt{-D'}$ has degree $D'$.
			\end{proof}
		\end{lemma}

		It is, however, possible that a CM elliptic curve $E/\bbQ$ may admit an $N$-congruence with a quadratic twist $E^d$ which is not isomorphic to $E^{-D}$. For later discussions (in particular when we give evidence for \thref{conj:allqtt}) it will be convenient to avoid CM elliptic curves since, for example, they do not satisfy Serre's open image theorem.

		We use the classification of images of the mod $p$ Galois representation for CM elliptic curves $E/\bbQ$ (see \cite[Propositions 1.14--1.16]{Z_OTPIOTMlRATECOQ}) and \thref{cong_q_twist} to classify nontrivial quadratic twist-type $p$-congruences between elliptic curves with CM (this is similar to the arguments of Cremona--Freitas \cite[Theorem 2.4 and Corollary 2.5]{CF_GMFTSTOCBEC}). We will also need to consider $2^k$-congruences. The image of the adelic Galois representation attached to a CM elliptic curve $E/\bbQ$ has been considered in \cite{LR_GRATECWCM} and \cite[Section 5]{DLRM_TACOEOGRATEC} and it would be interesting to complete the proof using these results. We find it easiest to argue using the curves $X_E^r(2^k)$. 
		
		We recall some basic facts about the curve $X_E^r(N)$. Let $E/K$ be an elliptic curve defined over a field of characteristic $0$. Then there exists a curve, $X_E^r(N)/K$, parametrising elliptic curves $(N, r)$-congruent to $E$ (see e.g., \cite[Section 1]{RS_FOECWCMPR}, \cite[Sections 4.4 and 4.5]{PSS_TOX7APSTX2Y3Z7}). More precisely there is a bijective correspondence between noncuspidal points $P \in X_E^r(N)(K) \setminus \{\text{cusps}\}$ and ($K$-isomorphism classes of) pairs $(E', \phi)$ where $\phi \colon E[N] \to E'[N]$ is an $(N, r)$-congruence. Note that by taking the $j$-invariant of $E'$ we obtain a natural map $j \colon X_E^r(N) \to X(1)$. The curve $X_E^r(N)$ is a twist of (any) geometrically connected component of the modular curve $X(N)$ as defined in Section \ref{sec:basic_defs}.

		\begin{table}
			\centering
			\begin{tabular}{c|c|c|c}
			$-D$ 			& Weierstrass equation for $E$		& $d$	& $(N, r)$  \\ 
			\hline
			$-3$			& $y^2 = x^3 + 1$					& $-1$	& $(4, 1)$ 	\\
							&									& $3$	& $(4, 3)$ 	\\
							& $y^2 = x^3 + 10$					& $5$	& $(10, 1)$	\\
							&									& $-15$	& $(10, 2)$	\\
							& $y^2 = x^3 + 98$					& $21$	& $(14, 1)$	\\
							&									& $-7$	& $(14, 3)$ \\
			\hline
			$-4$			& $y^2 = x^3 + ax$					& $a$	& $(4, 3)$	\\
							& $y^2 = x^3 + 2x$ 					& $2$	& $(8, 3)$	\\
							& $y^2 = x^3 - 2x$					& $2$	& $(8, 7)$	\\
			\hline
			$-7$			& $y^2 + xy = x^3 - x^2 - 2x - 1$	& $-1$	& $(4, 1)$	\\
							&									& $7$	& $(4, 3)$	\\
			\hline
			$-8$			& $y^2 = x^3 + 4x^2 + 2x$			& $-3$	& $(6, 1)$	\\
							&									& $6$	& $(6, 5)$	\\
							&									& $2$	& $(4, 3)$	\\
			\hline
			$-11$			& $y^2 + y = x^3 - x - 7x + 10$		& $-3$	& $(6, 1)$	\\
							&									& $33$	& $(6, 5)$	\\
			\hline
			$-12$			& $y^2 = x^3 - 15x + 22$			& $-1$	& $(4, 1)$	\\
							&									& $3$	& $(4, 3)$	\\
			\hline
			$-16$			& $y^2 = x^3 - 11x + 14$			& $2$	& $(4, 3)$	\\
			\hline
			$-28$			& $y^2 = x^3 - 595x + 5586$			& $-1$	& $(4, 1)$	\\
							&									& $7$	& $(4, 3)$	\\

			\end{tabular}
	
			\caption{All pairs $(E, d)$ where $E/\bbQ$ is an elliptic curve (defined up to quadratic twist) with CM by an imaginary quadratic order of discriminant $-D$ such that $E$ is nontrivially $(N, r)$-congruent with its quadratic twist by $d$.} \label{table:nontrivialCMcongs}
		\end{table}

		\begin{prop} \thlabel{prop:nontrivialCMcongs}
			Let $N > 2$ be an integer and let $E/\bbQ$ be an elliptic curve with CM by an imaginary quadratic order of discriminant $-D$.
			
			Let $d \in \bbZ$ be squarefree and suppose that $d \neq 1$ and $d \neq -1$ if $j(E) = 1728$ (i.e., suppose that $E$ is not isomorphic to $E^d$ over $\bbQ$). Then $E$ is $(N,r)$-congruent to $E^d$ if and only if either
			\begin{enumerate}[label=\eniii]
				\item $d$ is equal to $-D$ up to a square factor (i.e., the congruence is trivial), or
				\item $E$ is isomorphic to a quadratic twist of one of the elliptic curves appearing in Table \ref{table:nontrivialCMcongs} and $d$ and $(N,r)$ appear in the corresponding row.
			\end{enumerate}
			In particular no CM elliptic curve $E/\bbQ$ admits a nontrivial quadratic twist-type $(N, r)$-congruence when $N > 8$ and $N \neq 10, 14$.

			\begin{proof}
				It will suffice to consider the cases when $N = p$ is an odd prime and $N = 2^k$ for $k \geq 2$.

				Let $p$ be an odd prime. By \thref{cartan_cong_pl} and \cite[Theorem XI.2.3]{L_ITMF} if $E$ is $p$-congruent to a nontrivial quadratic twist $E^d$, then the image, $\bbP \bar{\rho}_{E,p}(G_{\bbQ})$, of $\bar{\rho}_{E,p}(G_{\bbQ})$ in $\PGL_2(\bbZ/p\bbZ)$ is dihedral and $\bbP \bar{\rho}_{E,p}(G_{\bbQ(\sqrt{d})})$ is a cyclic subgroup of index $2$. In particular if $\bbP \bar{\rho}_{E,p}(G_{\bbQ})$ is dihedral of order $> 4$ then it contains a unique cyclic subgroup of index $2$, in which case $E$ is $p$-congruent to a unique quadratic twist.

				Suppose that $j(E) \neq 0$. By \cite[Proposition 1.14]{Z_OTPIOTMlRATECOQ} for each odd prime $p$ not dividing $D$ we have that $\bar{\rho}_{E,p}(G_{\bbQ})$ is \emph{equal} to the normaliser of a Cartan subgroup of $\GL_2(\bbZ/p\bbZ)$ (which is split if and only if $-D$ is a quadratic residue modulo $p$) and $\bar{\rho}_{E,p}(G_{\bbQ(\sqrt{-D})})$ is the corresponding Cartan subgroup. Therefore unless $p = 3$ and $-D \equiv 1 \pmod{3}$ the image of $\bar{\rho}_{E,p}(G_{\bbQ})$ in $\PGL_2(\bbZ/p\bbZ)$ is dihedral of order $> 4$ hence $E$ is $p$-congruent to exactly one quadratic twist, namely $E^{-D}$, except when $D = 8$ or $11$ and $p = 3$. 
				
				When $D = 8$ or $11$ the image of the mod $3$ Galois representation, $\bar{\rho}_{E,3}(G_{\bbQ})$, is equal to the normaliser of a split Cartan subgroup and has exactly $3$ distinct index $2$ subgroups satisfying \thref{cong_q_twist}. In particular $E$ is nontrivially $(3,1)$-congruent to its quadratic twist by $-3$ and $(3, 2)$-congruent to its quadratic twist by $3D$ (and no other quadratic twists). For both $D = 8$ and $11$ the traces of Frobenius at $\ell = 17$ of $E$, $E^{-3}$, and $E^{-3D}$ are not congruent modulo $9$. Since the traces of Frobenius of $N$-congruent elliptic curves at good primes are congruent modulo $N$, none of these $3$-congruences are $9$-congruences. 

				If $p$ divides $D$ then by \cite[Proposition 1.14(iii)]{Z_OTPIOTMlRATECOQ} $\bar{\rho}_{E,p}(G_{\bbQ})$ in not contained in the normaliser of a Cartan subgroup, so by \thref{cartan_cong_pl} $E$ is not $p$-congruent to a nontrivial quadratic twist.

				Now suppose that $j(E) = 0$. Let $p \geq 5$ be a prime number. If $p = 5$ suppose that $E$ is not isomorphic to a quadratic twist of the elliptic curve $y^2 = x^3 + 10$ and if $p = 7$ suppose that $E$ is not isomorphic to a quadratic twist of the elliptic curve $y^2 = x^3 + 98$. It follows from \cite[Proposition 1.16(i)--(iv)]{Z_OTPIOTMlRATECOQ} that the the image of $\bar{\rho}_{E,p}(G_{\bbQ})$ in $\PGL_2(\bbZ/p\bbZ)$ is dihedral of order $> 4$ hence $E$ is $p$-congruent to a unique quadratic twist -- namely its $3$-isogenous quadratic twist $E^{-3}$.

				In the remaining two cases the image of the mod $p$ Galois representation, $\bar{\rho}_{E,p}(G_{\bbQ})$, (described in \cite[Proposition 1.16(ii)--(iii)]{Z_OTPIOTMlRATECOQ}) contains exactly three index $2$ subgroups, $H$, satisfying the conditions of \thref{cong_q_twist}. In particular the elliptic curve with Weierstrass equation $y^2 = x^3 + 10$ is $(5, 1)$-congruent with its quadratic twist by $5$ and $(5, 2)$-congruent with its quadratic twist by $-15$. Similarly the elliptic curve with Weierstrass equation $y^2 = x^3 + 98$ is $(7, 3)$-congruent with its quadratic twist by $-7$ and $(7, 1)$-congruent with its quadratic twist by $21$. In each of these cases by comparing traces of Frobenius at $\ell = 13$ we see that these $p$-congruences are not $p^2$-congruences.

				By \cite[Proposition 1.16(v)]{Z_OTPIOTMlRATECOQ} and \thref{cartan_cong_pl} if an elliptic curve $E$ with $j$-invariant $0$ is $3$-congruent to a nontrivial quadratic twist, then $\bar{\rho}_{E,p}(G_{\bbQ})$ is contained in a split Cartan subgroup of $\GL_2(\bbZ/3\bbZ)$ and $E$ is isomorphic to a quadratic twist of the elliptic curve with Weierstrass equation $y^2 = x^3 + 2$. In this case $E$ is both $(3, 1)$ and $(3, 2)$-congruent with its $3$-isogenous quadratic twist $E^{-3}$.

				We now determine the nontrivial quadratic twist-type $2^k$-congruences between CM elliptic curves $E/\bbQ$ where $k \geq 2$ using the curves $X_E^r(2^k)$. Silverberg \cite{S_EFOECWPMNR} and Fisher \cite{F_THOAG1C} computed equations for $X_E^r(4)$ and Z. Chen \cite{C_FOECWTSM8R} computed equations for $X_E^r(8)$.

				When $j(E) \not\in \{0, 1728\}$ we may take simultaneous quadratic twists of $E$ and $E^d$ so that $E$ is given by the corresponding Weierstrass equation in \cite[A.3]{S_ATITAOEC}. Similarly when $j(E) = 1728$ (i.e., $-D = -4$) we may assume that $E$ is given by the Weierstrass equation $y^2 = x^3 + ax$ where $a \in \bbZ$ is squarefree, and when $j(E) = 0$ (i.e., $-D = -3$) we may assume that $E$ is given by the Weierstrass equation $y^2 = x^3 + a$ where $a \in \bbZ$ is cubefree and $a > 0$. 
	
				In each case, by determining the rational points in the fibres of the maps $X_E^r(4) \to X(1)$ above $j(E)$ we see that the only nontrivial quadratic twist-type $(4, r)$-congruences between CM elliptic curves are those appearing in Table \ref{table:nontrivialCMcongs}.	For each such rational point $P \in X_E^r(4)(\bbQ)$ giving rise to a nontrivial quadratic twist-type $(4, r)$-congruence we compute the rational points in the fibre of the forgetful map $X_E^{r'}(8) \to X_E^r(4)$ above $P$ for each $r' \equiv r \pmod{4}$. We conclude that the only nontrivial quadratic twist-type $(8,r')$-congruences between CM elliptic curves are those in Table \ref{table:nontrivialCMcongs} (i.e., where $E$ is a quadratic twist of one of the elliptic curves $y^2 = x^3 \pm 2x$). In both of these cases comparing traces of Frobenius at $\ell = 29$ shows that the congruence is not a $16$-congruence.

				We illustrate this calculation in the case where $j(E) = 0$ and $r = 1$. The other cases are verified in \texttt{Magma} (for details see \cite{ME_ELECTRONIC}). Let $E$ be given by the Weierstrass equation $y^2 = x^3 + a$ where $a \in \bbZ$ is cubefree and $a > 0$. From Fisher's equations for $X_E^1(4) \cong \mathbb{P}^1$ \cite{F_THOAG1C} an elliptic curve is $(4, 1)$-congruent to $E$ if and only if it is isomorphic to an elliptic curve with Weierstrass equation
				\begin{equation*}
					y^2 = x^3 + 36 ast(27s^3 - 8 a t^3)(27s^3 + at^3) x + a (729 s^6 + 8 a^2 t^6) (729 s^6 - 2376 a s^3 t^3 - 8 a^2 t^6)
				\end{equation*}					
				for some point $(s : t) \in \mathbb{P}^1(\bbQ)$.
				
				If $a \neq 1$ the only rational points in the fibre of the $j$-map $X_E^1(4) \to X(1)$ above $0$ (i.e., where $st(27s^3 - 8 a t^3)(27s^3 + at^3) = 0$) are the points $(1:0)$ and $(0:1)$. The congruent elliptic curves corresponding to these points are $E$ itself and the elliptic curve with Weierstrass equation $y^2 = x^3 - a^5$ respectively. Note that the elliptic curve given by the Weierstrass equation $y^2 = x^3 - a^5$ is not a quadratic twist of $E$ when $a \neq 1$. Hence $E$ does not admit any quadratic twist-type $(4, 1)$-congruences.  
				
				Now suppose that $a = 1$. The rational points in the fibre of the map $X_E^1(4) \to X(1)$ above $0$ are the points $(1:0)$, $(-1 : 3)$, $(0:1)$, and $(2:3)$. The first two points correspond to $E$ itself and the last two correspond to the quadratic twist of $E$ by $-1$. 
				
				The elliptic curve $E : y^2 = x^3 + 1$ is $(8, 1)$-congruent to its quadratic twist by $-1$ if and only if there exists a rational point on the curve $X_E^1(8)$ which maps to either $(0:1)$ or $(2:3)$ on $X_E^1(4)$. By \cite[Theorem~1.1]{C_FOECWTSM8R} the curve $X_E^1(8)$ is isomorphic over $\bbQ$ to the curve in $\bbP^4$ given by the equations 
				\begin{align*}
					&2 x_1 x_3 + x_2^2 + 2 x_4^2 = 0,\\
					&x_3^2 - 2 x_1 x_2 - 2 x_0 x_4 = 0,\\
					&2 x_2 x_3 - x_1^2 + x_0^2 = 0,
				\end{align*}
				and with the forgetful map $X_E^1(8) \to X_E^1(4)$ given by $(x_0 : x_1 : x_2 : x_3 : x_4) \mapsto (x_0 : 3x_4)$. From these equations it is simple to check that there exist no rational points in the fibre of the map $X_E^1(8) \to X_E^1(4)$ above the points $(0:1)$ and $(2:3)$.
				
				Therefore for each $k \geq 2$ the only nontrivial quadratic twist-type $(2^k,1)$-congruences between elliptic curves of $j$-invariant $0$ are those recorded in Table \ref{table:nontrivialCMcongs}.

				In the general case, write $N = p_1^{l_1} ... p_n^{l_n}$ as a product of prime powers, where the primes $p_i$ are distinct. The claim follows immediately since an $N$-congruence between $E$ and a quadratic twist $E^d$ gives rise to a $p_i^{l_i}$-congruence between $E$ and $E^d$ for each $i = 1,...,n$. 
			\end{proof}
		\end{prop}

		\begin{remark} \thlabel{remark:weirdtrivialcongs}
			We proved in \thref{prop:nontrivialCMcongs} that the elliptic curve $E$ with Weierstrass equation $y^2 = x^3 + 2$ is $3$-congruent and $3$-isogenous to $E^{-3}$. This congruence does not arise from the restriction of the isogeny to the $3$-torsion. By considering the $p$-adic images of the Galois representation attached to CM elliptic curves, it may be possible to classify the $N$-congruences between a (CM) elliptic curve and its isogenous quadratic twist which do not arise from the restriction of an isogeny. We do not pursue this.
		\end{remark}

		Combining \thref{prop:nontrivialCMcongs} with the results of Section \ref{sec:technical} we deduce the case of \thref{thm:infmany} when $N$ is a prime power. For an odd prime $p$ let $\xi$ be a quadratic nonresidue modulo $p$.

		\begin{prop} \thlabel{thm:infmanypl}
			Let $p$ be a prime number and $l \geq 1$ be an integer. Let $\mathcal{S}$ be the set defined in \thref{thm:infmany}. Then there are infinitely many $j$-invariants of elliptic curves $E/\bbQ$ admitting nontrivial quadratic twist-type $(p^l,r)$-congruences (hence $(2p^l, r)$-congruences if $p$ is odd) if and only if $(p^l, r) \in \mathcal{S}$. 
			
			Furthermore, suppose that $(p^l, r) \not\in \mathcal{S}$ and $E/\bbQ$ is an elliptic curve admitting a nontrivial quadratic twist-type $(p^l, r)$-congruence. Then either $E$ is a quadratic twist of $E_{(32,3)a}$ or $E_{(32,3)b}$ (see \thref{thm:examples}) and $(p^l, r) = (32, 3)$, or $p^l \geq 19$ is odd, $r = -\xi$ and $j(E)$ is an exceptional $j$-invariant for $X(\ns {p^l}^+)$.

			\begin{proof}
				When $p$ is odd and $p^l \geq 13$ the modular curve $X(\ns {p^l}^+)$ has genus $\geq 2$ (see \cite[Section~5]{B_NONSCSMCATCNP}). Hence by Faltings' theorem and \thref{cartan_cong_pl}, when $p^l \geq 13$ there are at most finitely many $j$-invariants of elliptic curves $E/\bbQ$ admitting a quadratic twist-type $(p^l,-\xi)$-congruence. 
				
				By \cite[Theorem 1.1]{BPR_RPOX+0pr}, \cite{RSZB_LAIOGFECOQ}, and \cite{BDMSTV_ECKFTSCMCOL13} there are no exceptional points on $X(\s {p^l}^+)$ for $p^l \geq 9$. Moreover it is shown in \cite{BDMSTV_ECKFTSCMCOL13} and \cite{BDMTV_QCFMCAAE} that $X(\ns 13^+)$ and $X(\ns 17^+)$ have no exceptional points. Therefore by Propositions \ref{cartan_cong_pl} and \ref{prop:nontrivialCMcongs} there are no elliptic curves $E/\bbQ$ admitting a nontrivial quadratic twist-type $(9,2)$, $(11,2)$, $(13, 1)$, $(13, 2)$, $(17, 1)$, $(17, 3)$, or $(p^l, -1)$-congruence where $p^l \geq 19$. Moreover, if $E/\bbQ$ admits a nontrivial quadratic twist-type $(p^l, -\xi)$-congruence with $p^l \geq 19$ then there exists an exceptional point in the fibre of the $j$-map $X(\ns {p^l}^+) \to X(1)$ above $j(E)$.

				Similarly suppose that $p=2$ and $(p^l, r) \not\in \mathcal{S}$. Recall that the only exceptional points on $X(\ns 16^+)$ have $j$-invariants equal to $j(E_{(32,3)a})$ or $j(E_{(32,3)b})$ (this is due to Rouse and Zureick-Brown \cite{RZB_ECOQA2AIOG}, see the discussion in \thref{ex:32}). By Propositions \ref{prop:2adiccongs} and \ref{prop:nontrivialCMcongs} if $E/\bbQ$ admits a nontrivial quadratic twist-type $(p^l, r)$-congruence then $E$ is a quadratic twist of either $E_{(32,3)a}$ or $E_{(32,3)b}$ and $(p^l, r) = (32, 3)$. 

				Conversely suppose that $(p^l, r) \in \mathcal{S}$. By Propositions \ref{cartan_cong_pl} and \ref{prop:2adiccongs} there exists a subgroup $H^+ \subset \GL_2(\bbZ/p^l\bbZ)$ and $d \in \bbQ(X(H^+))$ such that for each rational point $P \in X(H^+)$ (except possibly for the finitely many zeroes and poles of $d$ and points above $j = \infty, 0, 1728$) an elliptic $E_P/\bbQ$ with $j$-invariant $j(P)$ is $(p^l, r)$-congruent to its quadratic twist by $d(P) \in \bbQ^\times$. 
				
				Note that exactly one of the following holds:
				\begin{enumerate}[label=\eniii]
					\item the fibre of the degree $2$ map $X(H) \to X(H^+)$ above $P$ has a $\bbQ$-rational point (equivalently $d(P)$ is a square),
					\item the curve $E_P$ is isogenous over $\bbQ$ to $E_P^{d(P)}$,
					\item the curve $E_P$ admits a nontrivial quadratic twist-type $(p^l, r)$-congruence with $E_P^{d(P)}$.
				\end{enumerate}	

				When $p = 2$ Rouse and Zureick-Brown \cite{RZB_ECOQA2AIOG} showed that each of the modular curves $X(H^+)$ (i.e., those corresponding to subgroups listed in Table \ref{table:2adiccongs}, except $X(\ns 16^+)$) are isomorphic over $\bbQ$ to $\bbP^1$. When $p$ is odd each of the modular curves $X(H^+)$ are isomorphic over $\bbQ$ to $\bbP^1$, except $X(\ns 11^+)$ which is an elliptic curve of rank $1$ (see e.g., \cite{SZ_MCOPPLWIMRP}). 

				When the curve $X(H^+)$ is isomorphic to $\mathbb{P}^1$ Hilbert's Irreducibility Theorem (see \cite[Chapter 9]{S_LOTMWT}) implies that there are infinitely many rational points $P \in X(H^+)(\bbQ)$ such that case (i) does not occur. The modular curve $X(\ns 11)$ has genus $4$ -- so case (i) occurs for at most finitely many rational points $P \in X(\ns 11^+)(\bbQ)$ by Faltings' Theorem.

				By \thref{lemma:CM_congs} case (ii) may only occur if $E_P$ has CM. But there are only finitely many $j$-invariants of elliptic curves over $\mathbb{Q}$ with CM, so (ii) occurs for finitely many rational points $P \in X(H^+)$.

				Thus (iii) holds for infinitely many rational points $P \in X(H^+)$. Since $X(H^+) \to X(1)$ has finite degree these correspond to infinitely many $j$-invariants, as required.
			\end{proof}
		\end{prop}

	\section{Double Covers of Modular Curves and Composite Level Quadratic Twist-Type Congruences} \label{sec:doublecovers}
	
		Let $p_1$ and $p_2$ be distinct odd primes and let $l_1,l_2 \geq 1$ be integers. We wish to construct elliptic curves admitting quadratic twist-type $p_1^{l_1} p_2^{l_2}$-congruences over $\bbQ$.
		
		It is natural to consider the modular curve $X(H_1^+, H_2^+)$, where $H_1^+ \subset \GL_2(\bbZ/p_1^{l_1}\bbZ) $ and $H_2^+ \subset \GL_2(\bbZ/p_2^{l_2}\bbZ)$ are the normalisers of Cartan subgroups. A noncuspidal point $P \in X(H_1^+, H_2^+)(K)$ gives rise to an elliptic curve $E/K$ admitting quadratic twist-type $p_1^{l_1}$ and $p_2^{l_2}$-congruences by \thref{cartan_cong_pl}.  However these congruences may be with non-isomorphic quadratic twists. We therefore consider a double cover of $X(H_1^+, H_2^+)$, which corresponds to requiring that the quadratic twists are indeed isomorphic. We make this formal in the following way:

		Let $N_1,N_2 > 2$ be coprime integers. Let $H_1^+ \subset \GL_2(\bbZ/N_1\bbZ)$ and $H_2^+ \subset \GL_2(\bbZ/N_2\bbZ)$ be subgroups, and suppose that there exist index $2$ subgroups $H_i \subset H_i^+$ and $g_i \in \GL_2(\bbZ/N_i\bbZ)$ satisfying the conditions of \thref{cong_q_twist}.

		Let $\iota_1$ and $\iota_2$ be the involutions giving the quotients $X(H_1) \to X(H_1^+)$ and $X(H_2) \to X(H_2^+)$. Note that while the involutions $\iota_i$ depend upon the groups $H_i^+$ containing $H_i$, in all the cases we consider, the subgroup $H^+$ is implicit in the notation for $H$. More precisely, if $p$ is a prime number and $H$ is a Cartan subgroup of $\GL_2(\bbZ/p^l\bbZ)$, then $H^+$ will be its normaliser. When $H^+ \subset \GL_2(\bbZ/2^k\bbZ)$ is one of the subgroups listed in Table \ref{table:2adiccongs} we write $H = H^+(\sqrt{d})$.
		
		Let $J_1 \colon X(H_1^+) \to X(1)$ and $J_2 \colon X(H_2^+) \to X(1)$ be the $j$-invariant maps. Let $\mathcal{E}_1/\bbQ(X(H_1^+))$ and $\mathcal{E}_2/\bbQ(X(H_2^+))$ be elliptic curves with $j$-invariants $J_1$ and $J_2$ respectively. By \thref{cong_q_twist}, there exist rational functions $d_1 \in \bbQ(X(H_1^+))^\times$ and $d_2 \in \bbQ(X(H_2^+))^\times$ (depending on the index $2$ subgroup $H_i$) such that $\mathcal{E}_i$ is $N_i$-congruent to its quadratic twist $\mathcal{E}_i^{d_i}$. 

		\begin{remark} \thlabel{rem:howtogetd}
			Let $F = \bbQ(X(H^+))$. For the computations in Sections \ref{sec:2kpl} and \ref{sec:bigconj} it will be necessary to explicitly compute these rational functions $d \in F^\times$ when $H^+ \subset \GL_2(\bbZ/p\bbZ)$ is the normaliser of a Cartan subgroup and $p$ is an odd prime number. 

			Suppose that $E/K$ is an elliptic curve admitting $(p, r)$-congruence with its quadratic twist by $d \in K^\times$. Then for $p \in \{3,5,7,11\}$ we can find $d$ from the equations for $X_E^r(p)$ by computing $K$-rational points on the fibre of the $j$-map $X_E^r(p) \to X(1)$ above $j(E)$ (equations for $X_E^r(p)$ and its $j$-invariant maps may be found in \cite{RS_FOECWCMPR} and \cite{F_THOAG1C} when $p = 3,5$ and in \cite{HK_SLCMXE7}, \cite{PSS_TOX7APSTX2Y3Z7}, and \cite{F_OFO7A11CEC} when $p = 7$ and \cite{F_OFO7A11CEC} when $p = 11$).

			However, since the genus of $X_E^r(p)$ grows quickly with $p$ it will instead be more convenient to note that by \thref{cartan_cong_pl}, $d \in F^\times/(F^\times)^2$ is the unique element such that the degree $2$ extension $\bbQ(X(H))/F$ is given by adjoining $\sqrt{d}$ to $F$. 
		\end{remark}
		
		We have the natural projections from $X(H_1^+, H_2^+)$, which give the commutative diagram 
		\[\begin{tikzcd}
			& {X(H_1^+, H_2^+)} \\
			X(H_1^+) && X(H_2^+) \\
			& {X(1)}
			\arrow["{\phi_1}"', from=1-2, to=2-1]
			\arrow["{\phi_2}", from=1-2, to=2-3]
			\arrow["{J_2}", from=2-3, to=3-2]
			\arrow["{J_1}"', from=2-1, to=3-2]
		\end{tikzcd}\]
	
		We define $\psi \colon \mathcal{C} \to X(H_1^+, H_2^+)$ to be the double cover such that the corresponding extension of function fields is given by adjoining $\sqrt{\phi_1^*(d_1) \phi_2^*(d_2)}$ to the function field $\bbQ(X(H_1^+, H_2^+))$ (note that $\mathcal{C}$ depends on the subgroups $H_1$ and $H_2$). Let $J \colon \mathcal{C} \to X(1)$ be the $j$-invariant map, and let $\mathcal{E}/\bbQ(\mathcal{C})$ be an elliptic curve with $j$-invariant $J$. Then $\mathcal{E}$ is $N_1$-congruent to its quadratic twist by $(\phi_1 \circ \psi)^*(d_1)$ and $N_2$-congruent to its quadratic twist by $(\phi_2 \circ \psi)^*(d_2)$. But $(\phi_1 \circ \psi)^*(d_1)$ and $(\phi_2 \circ \psi)^*(d_2)$ are equal up to a square factor, so $\mathcal{E}$ is $N_1N_2$-congruent to its quadratic twist by $(\phi_1 \circ \psi)^*(d_1)$.
		
		\begin{remark}
			Note that the curve $\mathcal{C}$ is modular. The involutions $\iota_1$ and $\iota_2$ induce quotients $X(H_1, H_2) \to X(H_1^+, H_2)$ and $X(H_1, H_2) \to X(H_1, H_2^+)$. Then the quotient $X(H_1, H_2) / (\iota_1 \iota_2)$ is isomorphic to $\mathcal{C}$. That is, we have a commutative diagram 
			\[\begin{tikzcd}
				& {X(H_1, H_2)} \\
				{X(H_1, H_2^+)} & {X(H_1, H_2)/\delta} & {X(H_1^+, H_2)} \\
				& {X(H_1^+, H_2^+)}
				\arrow["{\iota_1}"', from=1-2, to=2-1]
				\arrow["{\iota_2}", from=1-2, to=2-3]
				\arrow["\delta", from=1-2, to=2-2]
				\arrow[from=2-1, to=3-2]
				\arrow[from=2-3, to=3-2]
				\arrow[from=2-2, to=3-2]
			\end{tikzcd}\]
		\end{remark}
		
		Since the involutions $\iota_i$ will always be implicit in the notation of the subgroup $H_i$ we will drop the labels and denote the involution $\iota_1 \iota_2$ by $\delta$. We then have the following:
		
		\begin{lemma} \thlabel{prop:Xcong}
			Let $P \in X(H_1, H_2)/\delta $ be a non-cuspidal $K$-rational point. Let $\mathcal{E}_P$ be the specialisation of $\mathcal{E}$ at the point $P$ and suppose that $j(\mathcal{E}_P) \neq 0, 1728$. Then exactly one of the following holds:
			\begin{enumerate}[label=\eniii]
				\item the fibre of the degree $2$ map $X(H_1, H_2) \to X(H_1, H_2)/\delta$ above $P$ has a $K$-rational point,  
				\item the elliptic curve $\mathcal{E}_P$ admits a quadratic twist-type $N_1N_2$-congruence.
			\end{enumerate}
			
			Conversely, suppose that $E/K$ is $N_1N_2$-congruent to a nontrivial quadratic twist by $d \in K^\times$ and $j(E) \neq 0, 1728$. For each $i = 1, 2$ let $H_i^+ = \bar{\rho}_{E, N_i} (G_{K})$ and let $H_i = \bar{\rho}_{E, N_i} (G_{K(\sqrt{d})})$. Then there exists a $K$-rational point in the fibre of the map $X(H_1, H_2)/\delta \to X(1)$ above $j(E)$ which does not lift to a $K$-rational point on $X(H_1, H_2)$.

			\begin{proof}
				From \thref{lemma:modcurveimages} and the preceding discussion (and by specialising) we see that $\mathcal{E}_P$ is $N_1 N_2$-congruent to its quadratic twist by $(\phi_1 \circ \psi)^*(d_1)(P)$. Recall that by convention we say that $\mathcal{E}_P$ admits a quadratic twist-type $N_1N_2$-congruence if the quadratic twist is nontrivial. That is, if $(\phi_1 \circ \psi)^*(d_1)(P)$ is not a square in $K$. 
				
				The result then follows by noting that on the level of function fields the double cover $X(H_1, H_2) \to X(H_1, H_2)/\delta$ is given by adjoining the square root of $(\phi_1 \circ \psi)^*(d_1)$.

				The converse follows from \thref{lemma:modcurveimages} and the construction of the curve $\mathcal{C} \cong X(H_1, H_2)/\delta$, noting that the quadratic twists are equal up to a square factor.
			\end{proof}
		\end{lemma}

		\begin{coro} \thlabel{coro:classifyp1p2}
			Let $E/\bbQ$ be an elliptic curve with $j(E) \neq 0, 1728$. Let $p_1$ and $p_2$ be distinct odd primes. Write $N = p_1p_2$, and suppose that $r$ is an integer coprime to $N$. For each $i = 1,2$ let $H_i \subset \GL_2(\bbZ/p_i\bbZ)$ be a Cartan subgroup which is split if and only if $-r$ is a quadratic residue modulo $p_i$. 
			
			Then $E$ admits a quadratic twist-type $(N,r)$-congruence if and only if there exists a rational point in the fibre of the $j$-map $X(H_1, H_2)/\delta \to X(1)$ above $j(E)$ which does not lift to a rational point on $X(H_1, H_2)$. 
			
			Moreover the congruence is nontrivial if and only if $E$ does not have CM.

			\begin{proof}
				The first claim follows from \thref{cartan_cong_pl} and \thref{prop:Xcong}. The second claim follows from \thref{lemma:CM_congs} and \thref{prop:nontrivialCMcongs}.
			\end{proof}
		\end{coro}

		\begin{coro} \thlabel{coro:classify2kpl}
			Let $E/\bbQ$ be an elliptic curve without CM. Let $p$ be an odd prime, and let $l \geq 1$ and $k \geq 2$ be integers. Write $N = 2^k p^l$ and suppose that $r$ is an integer coprime to $N$. Let $H_2 \subset \GL_2(\bbZ/p^l\bbZ)$ be a Cartan subgroup which is split if and only if $-r$ is a quadratic residue modulo $p$.
			
			Then $E$ admits a quadratic twist-type $(N, r)$-congruence if and only if there exists a rational point in the fibre of the $j$-map $X(H_1, H_2)/\delta \to X(1)$ above $j(E)$ which does not lift to a rational point on $X(H_1, H_2)$, where $H_1 \subset \GL_2(\bbZ/2^k \bbZ)$ is one of the subgroups, $H$, listed in Table \ref{table:2adiccongs} which corresponds to the pair $(2^k, r)$. 
			
			\begin{proof}
				This follows from Propositions \ref{cartan_cong_pl} and \ref{prop:2adiccongs}, and \thref{prop:Xcong}.
			\end{proof}
		\end{coro}

	\section{Computing Modular Curves of Genus \texorpdfstring{$0$}{0} or \texorpdfstring{$1$}{1}} \label{sec:2kpl}
	
		Let $p$ be an odd prime, and let $k \geq 2$ and $l \geq 1$ be integers. By \thref{coro:classify2kpl} if an elliptic curve $E/\bbQ$ without CM admits a quadratic twist-type $2^kp^l$-congruence then $E$ gives rise to a rational point on a modular curve $X(H_1, H_2)/\delta$ where $H_1 \subset \GL_2(\bbZ/2^k\bbZ)$ is one of the subgroups, $H$, listed in Table \ref{table:2adiccongs} and $H_2 \subset \GL_2(\bbZ/p^l\bbZ)$ is a Cartan subgroup. 

		To complete the proof of \thref{thm:infmany} we compute explicit models for the modular curves $X(H_1, H_2)/\delta$ of genus $0$ or $1$. That is, those curves which (by Faltings' Theorem) could possibly give rise to infinite families of quadratic twist-type congruences over $\bbQ$.
		
		It is useful to note that if $H \subset \GL_2(\bbZ/p^l\bbZ)$ is a Cartan subgroup then for each $k \geq 1$ and we have the commutative diagram
		\[\begin{tikzcd}
			{X(\ns  {2^{k+1}}, H)/\delta} & {X(\ns  {2^{k+1}}^+)} \\
			{X(\ns  {2^{k}}, H)/\delta} & {X(\ns  {2^{k}}^+)} 
			\arrow[from=1-1, to=2-1]
			\arrow[from=1-1, to=1-2]
			\arrow[from=2-1, to=2-2]
			\arrow[from=1-2, to=2-2]
		\end{tikzcd}\]
		which realises $X(\ns  {2^{k+1}}, H)/\delta$ as a fibre product.

		\begin{remark}
			Since we compute equations by taking fibre products, we naturally obtain the $j$-maps, these may be accessed from the corresponding electronic data \cite{ME_ELECTRONIC}.
		\end{remark}

			\begin{lemma} \thlabel{curves:2kns3}
				The curves $X(\ns  2, \ns  3 )/\delta$ and $X(\ns  4, \ns  3 )/\delta$ are isomorphic over $\bbQ$ to $\bbP^1$.

				\begin{proof}
					We fix a model for $X(\ns  3^+) \cong \mathbb{P}^1$ with $j$-invariant map given by $t^3$ (see \cite[Th{\'e}or{\`e}me III.3.3]{H_CESLCM}). Let $\mathcal{E}/\bbQ(t)$ be an elliptic curve with $j$-invariant $t^3$. Now $\mathcal{E}$ is $(3,1)$-congruent to its quadratic twist by  
					\begin{equation} \label{eqn:dns3}
						-3(t^2 + 12t + 144),
					\end{equation}
					this may be seen from the discussion in \thref{rem:howtogetd} and the remark following Th{\'e}or{\`e}me III.3.5 in \cite[p. 505]{H_CESLCM}.

					Now $-3(t^2 + 12t + 144)\Delta(\mathcal{E}) \equiv -3(t-12)$ up to squares. Therefore $X(\ns  {2}, \ns  {3})/\delta$ is given by $s^2 = -3(t-12)$. After replacing $s$ by $6s$ we see that $X(\ns  {2}, \ns  {3})/\delta$ is isomorphic to $\bbP^1$ with $j$-invariant map
					\begin{equation} \label{eqn:jXns2ns3d}
						1728 \left(1 - s^2\right)^3.
					\end{equation}
					
					We fix a model (taken from \cite{RZB_ECOQA2AIOG}) for $X(\ns  4^+) \cong \bbP^1$ with $j$-invariant $(32 t_4 - 4)/t_4^4$. Recalling that $X(\ns  2^+) = X(1)$ the modular curve $X(\ns  {4}^+, \ns  {3}^+)/\delta$ is birational to the curve given by
					\begin{equation*}
						\frac{32 t_4 - 4}{t_4^4} = 1728(1 - s^2)^3.
					\end{equation*}
					Putting $t_4 = \frac{u}{6 (s^2 - 1)}$ gives a model for $X(\ns  {4}, \ns  {3})/\delta$ as 
					\begin{equation*}
						3s^2 = (u+1)^2(u^2 - 2u +3)
					\end{equation*}
					in $\mathbb{A}^2$. This curve has genus $0$, and we parametrise by setting
					\begin{equation*}
						(s, u) = \left( \frac{3 (t^2 - 2)(t^2 + 2) }{ (t^2 - 2t - 2)^2}, \frac{2 (t - 1) (t + 2)}{t^2 - 2t - 2} \right).
					\end{equation*}			
					This gives a model for $X(\ns  {4}, \ns  {3})/\delta \cong \bbP^1$ with forgetful map to $X(\ns  {4}^+)$ given by
					\begin{equation} \label{eqn:jns4ns3}
						t_4 = \frac{(t^2 - 2t - 2)^3}{24 (t^3 - 4)(t^3 + 2)}.
					\end{equation}	
				\end{proof}
			\end{lemma}

			\begin{lemma} \thlabel{curves:2ks3}
				We have the following:
				\begin{enumerate}[label=\eniii]
					\item The curve $X(\ns  2, \s  3 )/\delta$ is isomorphic over $\bbQ$ to $\bbP^1$.
					\item The curve $X(\ns  4, \s  3)/\delta$ has genus $1$ and is isomorphic over $\bbQ$ to the elliptic curve with LMFDB label \LMFDBLabel{48.a4}, which has rank $0$.
				\end{enumerate}
				Moreover the only rational points on $X(\ns  4, \s  3)/\delta$ are CM points which correspond to the discriminants $-3$ and $-11$.

				\begin{proof}
					We fix a model (see \cite[Th{\'e}or{\`e}me III.3.6(d)]{H_CESLCM}) for $X(\s  3^+) \cong \mathbb{P}^1$ with $j$-invariant map given by 
					\begin{equation} \label{eqn:js3}
						\left( \frac{3 t (t + 4)}{t + 3} \right)^3
					\end{equation}
					and let $\mathcal{E}/\bbQ(t)$ be an elliptic curve with this $j$-invariant. By \thref{rem:howtogetd} and \cite[Th{\'e}or{\`e}me III.3.6(b)]{H_CESLCM} we see that $\mathcal{E}$ is $(3,2)$-congruent to its quadratic twist by 
					\begin{equation}
						t^2 - 12.
					\end{equation}
		
					We have $(t^2 - 12)\Delta(\mathcal{E}) \equiv 3(t + 3)$ up to squares. Then by setting $3s^2 = t + 3$ we obtain an isomorphism over $\bbQ$ between $X(\ns  {2}, \s  {3})/\delta$ and $\bbP^1$ with $j$-invariant map $\left( \frac{3 (s^2 - 1) (3s^2 + 1)}{s^2} \right)^3$.
		
					Therefore $X(\ns  {4}, \s  {3})/\delta$ is birational to the curve given by
					\begin{equation*}
						\frac{32t_4 - 4}{t_4^4} = \left( \frac{3 (s^2 - 1) (3s^2 + 1)}{s^2} \right)^3.
					\end{equation*}
					Setting $t_4 = \frac{2 s^2}{3 (s^2 - 1) (3 s^2 + 1)}u$ we obtain a birational model for $X(\ns  {4}, \s  {3})/\delta$ defined by 
					\begin{equation*}
						4 s^2 u (u^3 - 4) = -3(s^2 - 1)(3s^2 + 1)
					\end{equation*}
					in $\mathbb{A}^2$. \texttt{Magma} verifies that 
					\begin{equation*}
						(s,u) = \left( \frac{6x^3 + x^2 y + 14 x^2 + 12 x y - 8 x + 8 y - 16}{(x + 4)^2 (2 x + y + 2)} , \frac{6(x^2 + 3 x + y + 2)}{(x + 4) (2 x + y + 2)} \right)
					\end{equation*}
					gives an explicit isomorphism from the elliptic curve with Weierstrass equation $y^2 = x^3 + x^2 - 4x - 4$ and LMFDB label \LMFDBLabel{48.a4} to $X(\ns  {4}, \s  {3})/\delta$.

					The Mordell-Weil group of the elliptic curve \LMFDBLabel{48.a4} is isomorphic to $(\bbZ/2\bbZ)^2$. From the $j$-maps above it is simple to check that each rational point on $X(\ns 4, \s 3)/\delta$ corresponds to an elliptic curve with CM by an order of discriminant $-3$ or $-11$.
				\end{proof}
			\end{lemma}

			\begin{lemma} \thlabel{curves:2kns5}
				The curve $X(\ns  2, \ns  5 )/\delta$ is isomorphic over $\bbQ$ to $\bbP^1$.

				\begin{proof}
					We fix a model (see \cite[Th{\'e}or{\`e}me III.5.13]{H_CESLCM}) for $X(\ns  5^+) \cong \mathbb{P}^1$ with $j$-invariant map given by 
					\begin{equation} \label{eqn:jns5}
						\frac{5^3 (t + 1) (2t + 1)^3 (2t^2 - 3t + 3)^3}{(t^2 + t - 1)^5}
					\end{equation}
					and let $\mathcal{E}/\bbQ(t)$ be an elliptic curve with this $j$-invariant. By \thref{rem:howtogetd} and the proof of Th{\'e}or{\`e}me III.5.18 in \cite[p. 568]{H_CESLCM} we see that $\mathcal{E}$ is $(5,2)$-congruent to its quadratic twist by 
					\begin{equation} \label{eqn:dns5}
						-(8t^2 - 12t + 7).
					\end{equation}
		
					We have $-(8t^2 - 12t + 7) \Delta(\mathcal{E}) \equiv -(t^2 + t - 1)$ up to squares. Therefore $X(\ns  {2}, \ns  {5})/\delta$ is given by $s^2 = -(t^2 + t - 1)$ in $\mathbb{A}^2$. We parametrise this conic by setting $(s,t) = (\frac{-u^2 + u + 1}{u^2 + 1}, \frac{-u^2 - 2u}{u^2 + 1})$ to obtain a model of $X(\ns  {2}, \ns  {5})/\delta \cong \bbP^1$ with $j$-invariant 
					\begin{equation}
						\frac{-5^3 (2u - 1) (4u^2 + u + 1)^3 (2u^2 + 3u + 3)^3 (u^2 + 4u - 1)^3}{(u^2 - u - 1)^5}.
					\end{equation}
				\end{proof}
			\end{lemma}

			\begin{lemma} \thlabel{curves:2ks5}
				The curve $X(\ns  2, \s  5 )/\delta$ is isomorphic over $\bbQ$ to the elliptic curve with LMFDB label \LMFDBLabel{20.a3}, which has rank $0$. Moreover the only noncuspidal rational points on $X(\ns  2, \s  5 )/\delta$ are CM points which correspond to the discriminants $-4$, $-11$, and $-19$.

				\begin{proof}
					We fix a model (see \cite[Th{\'e}or{\`e}me III.5.9]{H_CESLCM}) for $X(\s  5^+) \cong \mathbb{P}^1$ with $j$-invariant map given by 
					\begin{equation} \label{eqn:js5}
						\frac{(t + 3)^3(t^2 + t + 4)^3(t^2- 4t - 1)^3}{(t^2 + t - 1)^5}
					\end{equation}
					and let $\mathcal{E}/\bbQ(t)$ be an elliptic curve with this $j$-invariant. By \thref{rem:howtogetd} and the discussion in \cite[p. 548]{H_CESLCM} we see that $\mathcal{E}$ is $(5,1)$-congruent to its quadratic twist by 
					\begin{equation} \label{eqn:ds5}
						t^2 - 4t - 16.
					\end{equation} 
					
					We have $(t^2 - 4t - 16) \Delta(\mathcal{E}) \equiv t(t^2 + t -1)$ up to squares. Setting $x = t$ shows that $X(\ns  {2}, \s  {5})/\delta$ is isomorphic to the elliptic curve with Weierstrass equation $y^2 = x^3 + x^2 - x$ and LMFDB label \LMFDBLabel{20.a3}.

					The Mordell-Weil group of the elliptic curve \LMFDBLabel{20.a3} has order $6$. From the $j$-maps above it is simple to check that each point is either a cusp or a CM point corresponding to an elliptic curve with CM by an order of discriminant $-4$, $-11$, or $-19$.
				\end{proof}
			\end{lemma}

			\begin{lemma} \thlabel{curves:2k7}
				The curves $X(\ns  2, \ns  7)/\delta$ and $X(\ns  2, \s  7)/\delta$ are both isomorphic over $\bbQ$ to the elliptic curve with LMFDB label \LMFDBLabel{196.a2}, which has rank $1$.

				\begin{proof}
					We fix a model (see \cite[Th{\'e}or{\`e}me III.7.30]{H_CESLCM}) for $X(\ns  7^+) \cong \mathbb{P}^1$ with $j$-invariant map given by 
					\begin{equation} \label{eqn:jns7}
						\frac{(3t + 1)^3 (t^2 + 3t + 4)^3 (t^2 + 10t + 4)^3 (4t^2 + 5t + 2)^3}{(t^3 + t^2 - 2t - 1)^7}
					\end{equation}
					and let $\mathcal{E}/\bbQ(t)$ be an elliptic curve with this $j$-invariant. By \thref{rem:howtogetd} and \cite[Th{\'e}or{\`e}me III.7.32]{H_CESLCM} we see that $\mathcal{E}$ is $(7,1)$-congruent to its quadratic twist by 
					\begin{equation} \label{eqn:dns7}
						-(16t^4 + 68t^3 + 111t^2 + 62t + 11).
					\end{equation}
					
					We have $-(16t^4 + 68t^3 + 111t^2 + 62t + 11) \Delta(\mathcal{E}) \equiv -(t^3 + t^2 - 2t - 1)$ up to squares. Setting $x = -t$ shows that $X(\ns  {2}, \ns  {7})/\delta$ is isomorphic to the elliptic curve with Weierstrass equation $y^2 = x^3 - x^2 - 2x + 1$ and LMFDB label \LMFDBLabel{196.a2}.
		
					To compute a model for $X(\ns  {2}, \s  {7})/\delta$ we fix a model (see \cite[Th{\'e}or{\`e}me III.7.19]{H_CESLCM}) for $X(\s  7^+) \cong \mathbb{P}^1$ with $j$-invariant map given by 
					\begin{equation} \label{eqn:js7}
						\frac{(1-t) (t-2)^3 (t^2 + 3t -3)^3 (t^2 + 3t + 4)^3 (t^4 + t^3 - t^2 + 2t + 4)^3}{(t^3 + t^2 -2t - 1)^7}
					\end{equation}
					and let $\mathcal{E}/\bbQ(t)$ be an elliptic curve with this $j$-invariant. By \thref{rem:howtogetd} and \cite[Proposition III.7.24]{H_CESLCM} we see that $\mathcal{E}$ is $(7,6)$-congruent to its quadratic twist by 
					\begin{equation}
						(t^4 + 6t^3 + 3t^2 - 18t - 19).
					\end{equation}
					
					We have $(t^4 + 6t^3 + 3t^2 - 18t - 19) \Delta(\mathcal{E}) \equiv -(t^3 + t^2 - 2t - 1)$ up to squares. Setting $x = -t$ shows that $X(\ns  {2}, \s  {7})/\delta$ is isomorphic to the elliptic curve with Weierstrass equation $y^2 = x^3 - x^2 - 2x + 1$ and LMFDB label \LMFDBLabel{196.a2}.
				\end{proof}
			\end{lemma}

			\begin{lemma} \thlabel{curves:2kns9}
				The curve $X(\ns  2, \ns  9)/\delta$ is isomorphic over $\bbQ$ to the elliptic curve with LMFDB label \LMFDBLabel{324.a2}, which has rank $1$.

				\begin{proof}
					We fix a model (due to Baran \cite[Proposition 4.6]{B_AMCOL9ATCNOP}) for $X(\ns  9^+) \cong \mathbb{P}^1$ with forgetful map $X(\ns  9^+) \to X(\ns  3^+)$ given by
					\begin{equation}
						t_3 = \frac{-3(t^3 + 3t^2 - 6t + 4) (t^3 + 3t^2 + 3t + 4) (5t^3 - 3t^2 - 3t + 2)}{(t^3 - 3t + 1)^3}.
					\end{equation}
					Once again note that we may realise $X(\ns  {2}, \ns  {9})/\delta$ as the fibre product
					\[\begin{tikzcd}
						{X(\ns  {2}, \ns  {9})/\delta} & {X(\ns  9^+)} \\
						{X(\ns  {2}, \ns  {3})/\delta} & {X(\ns  3^+)} 
						\arrow[from=1-1, to=2-1]
						\arrow[from=1-1, to=1-2]
						\arrow[from=2-1, to=2-2]
						\arrow[from=1-2, to=2-2]
					\end{tikzcd}\]
					Noting (\ref{eqn:jXns2ns3d}) the modular curve $X(\ns  {2}, \ns  {9})/\delta$ is given by 
					\begin{equation*}
						12(1-s^2) = \frac{-3(t^3 + 3t^2 - 6t + 4) (t^3 + 3t^2 + 3t + 4) (5t^3 - 3t^2 - 3t + 2)}{(t^3 - 3t + 1)^3}.
					\end{equation*}
					After making the change of variables $s = \frac{3(t + 1)(t^2 - t + 1)(t^2 + 2t - 2)}{2 (t^3 - 3 t + 1)^2} u$ we have a model for $X(\ns  {2}, \ns  {9})/\delta$ as 
					\begin{equation*}
						(t + 1) u^2 = (t^3 - 3 t + 1)
					\end{equation*}
					in $\mathbb{A}^2$. Then setting $x = \frac{3}{t + 1}$ and $y = \frac{-3u}{t + 1}$ shows that $X(\ns  {2}, \ns  {9})/\delta$ is isomorphic to the elliptic curve with Weierstrass equation $y^2 = x^3 - 9x + 9$ and LMFDB label \LMFDBLabel{324.a2}. The forgetful map $X(\ns  {2}, \ns  {9})/\delta \to X(\ns  9^+)$ is given by $t_9 = \frac{3-x}{x}$.
				\end{proof}
			\end{lemma}

			\begin{lemma} \thlabel{curves:Gmisc}
				We have the following:
				\begin{enumerate}[label=\eniii]
					\item The curve $X(G_{11}(\sqrt{-1}), \ns 3)/\delta$ is isomorphic over $\bbQ$ to the elliptic curve with LMFDB label \LMFDBLabel{72.a4} and rank $0$. Moreover the only rational points on the modular curve $X(G_{11}(\sqrt{-1}), \ns 3)/\delta$ are CM points corresponding to the discriminant $-4$.
					
					\item The curve $X(G_{9}, \ns 3^+)$ is isomorphic over $\bbQ$ to the elliptic curve with LMFDB label \LMFDBLabel{144.a3} and rank $0$. Moreover the only rational points on $X(G_{9}, \ns 3^+)$ are a cusp and a CM point corresponding to the discriminant $-4$. In particular the curves $X(G_9(\sqrt{-1}), \ns 3)/\delta$ and $X(G_9(\sqrt{-1}), \s 3)/\delta$ have no exceptional points.
				\end{enumerate}

				\begin{proof}
					(i) To compute a model for $X(G_{11}(\sqrt{-1}), \ns  {3})/\delta$ we fix a model (computed by Rouse and Zureick-Brown \cite{RZB_ECOQA2AIOG}) for $X(G_{11}) \cong \mathbb{P}^1$ with $j$-invariant map given by $\frac{(t^2 - 16)^3}{t^2}$. We then see that $X(G_{11}, \ns  {3}^+)$ is given by requiring that $\frac{(t^2 - 16)^3}{t^2}$ is a cube. In particular setting $t = s^3$ gives an isomorphism between $X(G_{11}, \ns  {3}^+)$ and $\bbP^1$ with forgetful map to $X(\ns 3^+)$ given by 
					\begin{equation} \label{eqn:jg11ns3}
						t_3 = \frac{(s^6 - 16)}{s^2}.
					\end{equation}
					
					Now $3(t_3^2 + 12t_3 + 144) \equiv 3(s^2 - 2s + 4)(s^2 + 2s + 4)$ up to a square factor. By (\ref{eqn:dns3}) the modular curve $X(G_{11}(\sqrt{-1}), \ns  {3})/\delta$ is birational to the curve
					\begin{equation*}
						u^2 = 3(s^2 - 2s + 4)(s^2 + 2s + 4)
					\end{equation*} 
					in $\mathbb{A}^2$. Taking 
					\begin{equation*}
						(s,u) = \left( \frac{-2(3x + y + 15)}{3x - y + 15}, \frac{12(x - 1)(x + 5)(x + 11)}{(3x - y + 15)^2} \right)
					\end{equation*}
					gives an explicit isomorphism from the elliptic curve with Weierstrass equation $y^2 = x^3 - 39x - 70$ and LMFDB label \LMFDBLabel{72.a4} to $X(G_{11}(\sqrt{-1}), \ns  {3})/\delta$. The elliptic curve \LMFDBLabel{72.a4} has Mordell-Weil group isomorphic to $(\bbZ/2\bbZ)^2$. From the $j$-maps above it is simple to check that all the rational points on $X(G_{11}(\sqrt{-1}), \ns  {3})/\delta$ correspond to an elliptic curve with CM by an order of discriminant $-4$.

					(ii) To compute a model for $X(G_{9}, \ns  {3}^+)$ we fix a model for $X(G_{9}) \cong \mathbb{P}^1$ (computed by Rouse and Zureick-Brown \cite{RZB_ECOQA2AIOG}) with $j$-invariant map given by $\frac{(t^2 + 48)^3}{t^2 + 64}$. We then see that $X(G_{9}, \ns  {3}^+)$ is given by requiring that $\frac{(t^2 + 48)^3}{t^2 + 64}$ is a cube. Therefore $X(G_{9}, \ns  {3}^+)$ is birational to the curve given by
					\begin{equation*}
						t^2 + 64 = s^3
					\end{equation*}
					in $\mathbb{A}^2$. Setting $t = 8y$ and $s = 4x$ we see that $X(G_{9}, \ns  {3}^+)$ is isomorphic to the elliptic curve with Weierstrass equation $y^2 = x^3 - 1$ and LMFDB label \LMFDBLabel{144.a3}. The elliptic curve \LMFDBLabel{144.a3} has Mordell-Weil group isomorphic to $\bbZ/2\bbZ$. From the $j$-invariant maps above it is simple to check that one point on $X(G_{9}, \ns 3^+)$ is a cusp and the other corresponds to an elliptic curve with CM by an order of discriminant $-4$. 
					
					The curve $X(G_9(\sqrt{-1}), \ns 3)/\delta$ is a double cover of $X(G_{9}, \ns 3^+)$ and $X(G_9(\sqrt{-1}), \s 3)/\delta$ is a degree $4$ cover of $X(G_{9}, \ns 3^+)$ via the factorisation $X(G_9(\sqrt{-1}), \s 3)/\delta \to X(G_9(\sqrt{-1}), \s 3^+) \to X(G_{9}, \ns 3^+)$. In particular, the only rational points on the curves $X(G_9(\sqrt{-1}), \ns 3)/\delta$ and $X(G_9(\sqrt{-1}), \s 3)/\delta$ are CM points and cusps.
				\end{proof}
			\end{lemma}

			\begin{remark}
				Note that we have omitted calculating models for $X(G_9(\sqrt{-\Delta}), \ns 3)/\delta$, $X(G_9(\sqrt{-\Delta}), \s 3)/\delta$, $X(G_{11}(\sqrt{-\Delta}), \ns 3)/\delta$ and $X(G_{11}(\sqrt{-\Delta}), \s 3)/\delta$ because in \thref{curves:2kns3} and \thref{curves:2ks3} we showed that $X(\ns  2, \ns  3 )/\delta$ and $X(\ns  2, \s  3 )/\delta$ have infinitely many rational points. By \thref{coro:classify2kpl} each of these curves may only give rise to quadratic twist-type $(12, 7)$ or $(12, 11)$-congruences. 
			\end{remark}
			
			\begin{table}
				\centering
				\begin{tabular}{c|c|c|c|c|c}
				$(2^k p^l,r)$ 	& Modular Curve 										& Genus 	& LMFDB Label 	& Rank	& Lemma					\\ 
				\hline
				$(12,1)$		& $X(G_9(\sqrt{-1}), \ns 3)/\delta$						& $1$		& \LMFDBLabelIsog{144.a}	& $0$	& \ref{curves:Gmisc}		\\
								& $X(G_{11}(\sqrt{-1}), \ns 3)/\delta$					& $1$		& \LMFDBLabel{72.a4}		& $0$	& \ref{curves:Gmisc}		\\
				$(12,5)$		& $X(G_9(\sqrt{-1}), \s 3)/\delta$						& $1$		& \LMFDBLabelIsog{144.a}	& $0$	& \ref{curves:Gmisc}		\\
				$(12,7)$		& $X(\ns 2, \ns 3)/\delta$								& $0$		&					&		& \ref{curves:2kns3}		\\
				$(12,11)$		& $X(\ns 2, \s 3)/\delta$								& $0$		&					&		& \ref{curves:2ks3}		\\
				$(20,3)$		& $X(\ns 2, \ns 5)/\delta$								& $0$		&					&		& \ref{curves:2kns5}		\\
				$(20,11)$		& $X(\ns 2, \s 5)/\delta$								& $1$		& \LMFDBLabel{20.a3}		& $0$	& \ref{curves:2ks5}		\\
				$(24,11)$		& $X(\ns 4, \s 3)/\delta$								& $1$		& \LMFDBLabel{48.a4}		& $0$	& \ref{curves:2ks3}		\\
				$(24,19)$		& $X(\ns 4, \ns 3)/\delta$								& $0$		& 					&		& \ref{curves:2kns3} 		\\
				$(28,3)$		& $X(\ns 2, \s 7)/\delta$								& $1$		& \LMFDBLabel{196.a2}	& $1$	& \ref{curves:2k7}		\\
				$(28,11)$		& $X(\ns 2, \ns 7)/\delta$								& $1$		& \LMFDBLabel{196.a2}	& $1$	& \ref{curves:2k7}		\\
				$(36,19)$		& $X(\ns 2, \ns 9)/\delta$								& $1$		& \LMFDBLabel{324.a2}	& $1$	& \ref{curves:2kns9}
				\end{tabular}
		
				\caption{Modular curves of genus $\leq 1$ which parametrise quadratic twist-type $(2^kp^l,r)$-congruences of elliptic curves where $k \geq 2$, $l \geq 1$ and $p$ is an odd prime.} \label{table:infmanycurves}
			\end{table}

			With these explicit models we now conclude the proof of \thref{thm:infmany}.

			\begin{proof}[Proof of \thref{thm:infmany}]
				We have already proved the case where $N$ is a prime power in \thref{thm:infmanypl}. It therefore suffices to consider the case where $N$ is divisible by two or more distinct primes. When $N = 2p^l$ and $p$ is an odd prime, the claim follows immediately from \thref{thm:infmanypl} (an elliptic curve is $2$-congruent to all of its quadratic twists).
				
				First suppose that $(N, r)$ is not in the set $\mathcal{S}$ in the statement of \thref{thm:infmany}. Suppose that $N$ is divisible by distinct \emph{odd} prime numbers $p_1$ and $p_2$. If an elliptic curve $E/\bbQ$ with $j(E) \neq 0, 1728$ admits a quadratic twist-type $N$-congruence then by \thref{coro:classifyp1p2} there exists a rational point on one of the $4$ modular curves $X(H_1, H_2)/\delta$ above $j(E)$, where $H_i \subset \GL_2(\bbZ/p_i\bbZ)$ is a Cartan subgroup. Table \ref{table:pqgenera} lists the genera of the modular curves $X(H_1, H_2)/\delta$ whenever $X(H_i^+)$ has infinitely many rational points. It follows from Faltings' Theorem that for each such $N$ there are at most finitely $j$-invariants of elliptic curves $E/\bbQ$ admitting a quadratic twist-type $N$-congruence. 

				Now suppose that $N$ is divisible by $2^k p^l$ where $p$ is an odd prime and $l \geq 1$ and $k \geq 2$ are integers. Table \ref{table:2kgenera} records the genera of the modular curves $X(H_1, H_2)/\delta$ from \thref{coro:classify2kpl} whenever $X(H_i^+)$ has infinitely many rational points. It follows from Faltings' Theorem that if there exist infinitely many $j$-invariants of elliptic curves admitting quadratic twist-type $(N,r)$-congruences then $(N,r) \in \mathcal{S}$ or there are infinitely many points on one of the modular curves $X(G_{11}(\sqrt{-1}), \ns 3)/\delta$, $X(G_{9}(\sqrt{-1}), \ns 3)/\delta$, $X(G_{9}(\sqrt{-1}), \s 3)/\delta$, $X(\ns 2, \s 5)/\delta$,  or $X(\ns 4, \s 3)/\delta$. But each of these curves has finitely many rational points (see Table \ref{table:infmanycurves}).

				Conversely, as recorded in Table \ref{table:infmanycurves}, the modular curves $X(\ns  2, \ns  3 )/\delta$, $X(\ns  2, \s  3 )/\delta$, $X(\ns  4, \ns  3 )/\delta$, $X(\ns  2, \ns  5 )/\delta$, $X(\ns  2, \ns  7)/\delta$, $X(\ns  2, \s  7)/\delta$, and $X(\ns  2, \ns  9)/\delta$ have infinitely many rational points. 

				For a rational point $P$ on one of these modular curves (not lying above $j = \infty, 0, 1728$) let $E_P/\bbQ$ be an elliptic curve with $j$-invariant $j(P)$. Then exactly one of the following holds:
				\begin{enumerate}[label=\eniii]
					\item the fibre of the degree $2$ map $X(\ns {2^k}, H) \to X(\ns {2^k}, H)/\delta$ above $P$ has a $\bbQ$-rational point (equivalently $\Delta(E_P)$ is a square),
					\item the curve $E_P$ is isogenous over $\bbQ$ to $E_P^{\Delta(E_P)}$,
					\item the curve $E_P$ admits a nontrivial quadratic twist-type congruence with $E_P^{\Delta(E_P)}$.
				\end{enumerate}				
				In the cases where the modular curve $X(\ns {2^k}, H)/\delta$ has genus $0$, Hilbert's Irreducibility Theorem (see \cite[Chapter 9]{S_LOTMWT}) implies that for infinitely many rational points $P$ case (i) does not occur. In the cases $X(\ns  2, \ns  7)/\delta$, $X(\ns  2, \s  7)/\delta$, and $X(\ns  2, \ns  9)/\delta$ the corresponding double cover in (i) has genus $5$, $5$, and $7$ respectively - so case (i) occurs for at most finitely many rational points by Faltings' Theorem.

				By \thref{lemma:CM_congs} (ii) may only occur if $E_P$ has CM. But there are only finitely many $j$-invariants of elliptic curves over $\mathbb{Q}$ with CM, so (ii) occurs for at most finitely many points $P \in X(H^+)$. 
				
				Therefore (iii) occurs for infinitely many rational points $P \in X(\ns  2^k, H)/\delta$. Hence when $(N, r) \in \mathcal{S}$ it follows from \thref{coro:classify2kpl} that there exist infinitely many $j$-invariants of elliptic curves admitting nontrivial quadratic twist-type $(N,r)$-congruences over $\bbQ$, as required.
			\end{proof}
	
	\section{Evidence for \texorpdfstring{\thref{conj:allqtt}}{Conjecture 1.4}} \label{sec:bigconj}
	
		This section is motivated by two goals. First we wish to construct examples of nontrivial quadratic twist-type $(N,r)$-congruences for some $(N,r) \not\in \mathcal{S}$ (in particular we wish to prove \thref{thm:examples}, see Examples \ref{ex:30cong} and \ref{ex:48cong}). Second we provide evidence for \thref{conj:allqtt}.

		By \thref{thm:infmanypl}, the prime power case of \thref{conj:allqtt} is equivalent to the conjecture that there exist no exceptional points on the modular curve $X(\ns {p^l}^+)$ for all odd prime powers $p^l \geq 19$ (see e.g., \cite[Conjecture 1.5]{RSZB_LAIOGFECOQ}). This is our ``Assumption (A)'' in \thref{thm:conditional}.
		
		If we assume (A) then by Propositions \ref{prop:nontrivialCMcongs} and \ref{thm:infmanypl} when $N$ is not a prime power (or twice an odd prime power) and $(N,r) \not\in \mathcal{S}$, nontrivial quadratic twist-type $(N,r)$-congruences may only arise from exceptional points on the modular curves $X(H_1, H_2)/\delta$ in Corollaries \ref{coro:classifyp1p2} and \ref{coro:classify2kpl} where $X(H_i^+)$ has infinitely many rational points. These curves, and their genera, are listed in Tables \ref{table:pqgenera} and \ref{table:2kgenera}.
		
		We therefore compute explicit models for the modular curves $X(H_1, H_2)/\delta$ in Tables \ref{table:pqgenera} and \ref{table:2kgenera} of genus $\leq 4$ whose noncuspidal rational points give rise to elliptic curves admitting an $(N,r)$-congruence with a quadratic twist, where $(N,r) \not\in \mathcal{S}$. Each such curve is listed (together with their known exceptional $j$-invariants) in Table \ref{table:ratpoints}.

		This allows us to prove the following conditional form of \thref{conj:allqtt}:
				
		\begin{theorem} \thlabel{thm:conditional}
			Assume the following:
			\begin{enumerate}[label=\enABC]
				\item There are no exceptional points on the modular curves $X(\ns {p^l}^+)$ for all odd prime powers $p^l \geq 19$,
				\item There are no (non-CM) points of canonical height $\geq 10^{15}$ on the genus $1$ modular curves in Table \ref{table:models} lifting to rational points on the double cover $X(H_1, H_2)/\delta$,
				\item There are no exceptional points on the modular curves in Tables \ref{table:pqgenera} and \ref{table:2kgenera} of genus $\geq 5$.
			\end{enumerate}
		
			Then \thref{conj:allqtt} holds.
		\end{theorem}		

		\begin{table}
			\centering
				\begin{tabular}{c|c|c|c}
					$(N,r)$		& Modular Curve		 						& Exceptional $j$-Invariants			& Reference 			\\
					\hline
					$(12, 1)$	& $X(G_9(\sqrt{-1}), \ns 3)/\delta$			& 									& \thref{curves:Gmisc}		\\
								& $X(G_{11}(\sqrt{-1}), \ns 3)/\delta$		&									& \thref{curves:Gmisc}		\\
					$(12, 5)$	& $X(G_9(\sqrt{-1}), \s 3)/\delta$			&									& \thref{curves:Gmisc}		\\
								& $X(G_{11}(\sqrt{-1}), \s 3)/\delta$		&									& \thref{curves:Gmisc4}		\\
					$(20, 11)$	& $X(\ns 2, \s 5)/\delta$					&									& \thref{curves:2ks5}			\\
					$(24, 7)$	& $X(\s 4, \ns 3)/\delta$					&									& \thref{curves:s4ns3}		\\
					$(24,11)$	& $X(\ns 4, \s 3)/\delta$					&									& \thref{curves:2ks3}		\\
					$(30, 1)$	& $X(\ns  3, \s  5)/\delta$					&									& Assumption (B) 			\\
					$(30, 7)$	& $X(\ns  3, \ns  5)/\delta$				& $\frac{- 3^3 \cdot 5^6 \cdot 199^3 \cdot 809^3 \cdot 5059^3}{61^{15}}$
																												& \thref{prop:ns3ns5pts}	\\
					$(40, 3)$	& $X(\ns  4, \ns  5)/\delta$				&									& \thref{curves:ns4ns5}	\\
					$(42, 1)$	& $X(\ns  3, \ns  7)/\delta$				&									& Assumption (B)			\\
					$(44, 3)$	& $X(\ns  2, \ns  11)/\delta$				&									& Assumption (B)			\\
					$(48, 19)$	& $X(\ns  8, \ns  3)/\delta$				& $\frac{2^{15} \cdot 3^3 \cdot 5^3 \cdot 31^3 \cdot 41^3 \cdot 47^3 \cdot 83^3 \cdot 293^3}{23^{24}}$
																												& Assumption (B)			\\
				\end{tabular}
			\caption{Exceptional $j$-invariants for the modular curves $X(H_1, H_2)/\delta$ in Tables \ref{table:pqgenera} and \ref{table:2kgenera} with genus $\leq 4$ and whose exceptional points give rise to elliptic curves admitting an $(N,r)$-congruence with a quadratic twist, where $(N,r) \not\in \mathcal{S}$.} \label{table:ratpoints}
			
		\end{table}

		\begin{table}
			\centering
			\begin{adjustbox}{width=\textwidth}
				\begin{tabular}{c|c|c|c|c}
					Modular Curve		 								& Affine Equation in $\mathbb{A}^3$			& Genus $1$ Modular Curve						& LMFDB Label						& Lemma		\\ 
					\hline
					\multirow{2}{*}{$X(\ns  3, \ns  5)/\delta$}			& $v^2 + v = u^3 + 1$						& \multirow{2}{*}{$X(\ns  3^+, \ns 5^+)$}		&\multirow{2}{*}{\LMFDBLabel{225.c2}}	& \multirow{2}{*}{\ref{curves:ns3ns5}}	\\
					& $w^2 = 3(3u^2 - 2u - 4v + 3)(9u^2 - 10u + 5)$		&	& 	&		\\
					\hline
					\multirow{2}{*}{$X(\ns  3, \s  5)/\delta$}			& $v^2 + v = u^3 + 1$						& \multirow{2}{*}{$X(\ns  3^+, \s 5^+)$}		&\multirow{2}{*}{\LMFDBLabel{225.c2}}	& \multirow{2}{*}{\ref{curves:ns3s5}}\\
					& $w^2 = -(3v + 3u^3 - 12u^2 + 12u - 6)$	&	& 	&		\\
					\hline
					\multirow{2}{*}{$X(\ns  3, \ns  7)/\delta$}			& $v^2 + v = u^3 + 12$						& \multirow{2}{*}{$X(\ns  3^+, \ns 7^+)$}		&\multirow{2}{*}{\LMFDBLabel{441.d2}}	& \multirow{2}{*}{\ref{curves:ns3ns7} }	\\
					& $w^2 = 4u^3 + 21u^2 - 42u + 49$			&	& 	&		\\
					\hline
					\multirow{2}{*}{$X(\ns  2, \ns  11)/\delta$}		& $v^2 + v = u^3 - u^2 - 7u + 10$			& \multirow{2}{*}{$X(\ns  11^+)$}				&\multirow{2}{*}{\LMFDBLabel{121.b2}}	& \multirow{2}{*}{ \ref{curves:2kns11}}	\\
					& $w^2 = (-3u^2 + 24u - 37)v - u^4 + 5u^3 + 18u^2 - 95u + 94$	&	& 	&		\\
					\hline
					\multirow{2}{*}{$X(\ns  8, \ns  3)/\delta$}			& $v^2 = u^3 + 8$							& \multirow{2}{*}{$X(\ns  8^+, \ns 3^+)$}		&\multirow{2}{*}{\LMFDBLabel{576.e4}}	& \multirow{2}{*}{\ref{curves:ns8ns3}}	\\
					& $w^2 = 3(3u^3 + 3u^2 + 2uv - 8u + 8)$		&	& 	&		\\
				\end{tabular}
			\end{adjustbox}
			\caption{Affine models for the modular curves in Table \ref{table:ratpoints} which naturally arise as a double cover of a modular curve of genus $1$ and rank $1$. The double covers are given by $(u,v,w) \mapsto (u,v)$.} \label{table:models}
			
		\end{table}

		\begin{remark}
			Assumption (A) in \thref{thm:conditional} would follow from \cite[Conjecture 1.5]{RSZB_LAIOGFECOQ} together with \cite[Corollary 1.3]{BDMSTV_ECKFTSCMCOL13} and \cite[Theorem 1.4]{BDMTV_QCFMCAAE}. In \thref{prop:ns3ns5pts} we prove that the only rational points on the curve $X(\ns 3, \ns 5)/\delta$ are those listed in Table \ref{table:pts_ns3ns5d}, so assumption (B) may be weakened slightly to exclude $X(\ns 3, \ns 5)/\delta$.
		\end{remark}

		\subsection{Computing Modular Curves of Genus \texorpdfstring{$\leq 4$}{leq 4}} \label{sec:gle4}
		
			We first show that the modular curves in Table \ref{table:ratpoints} which do not appear in Table \ref{table:models} arise naturally as a double cover of a modular curve with no exceptional points. Note that the curves in Table \ref{table:ratpoints} of genus $\leq 1$ were computed in Section \ref{sec:2kpl}.

			\begin{lemma} \thlabel{curves:Gmisc4}
				The curve $X(G_{11}, \s 3^+)$ is isomorphic over $\bbQ$ to the elliptic curve with LMFDB label \LMFDBLabel{48.a4} and rank $0$. Moreover the modular curves $X(G_{11}, \s 3^+)$ and $X(G_{11}(\sqrt{-1}), \s 3)/\delta$ have no noncuspidal rational points.

				\begin{proof}
					Recalling (\ref{eqn:js3}) we see that the double cover $X(\s 3^+) \to X(\ns 3^+)$ is given on the level of function fields by setting $\frac{3 t (t + 4)}{t + 3} = t_3$. In particular by (\ref{eqn:jg11ns3}) the modular curve $X(G_{11}, \s 3^+)$ is birational to the curve given by 
					\begin{equation*}
						\frac{3 t (t + 4)}{t + 3} = \frac{(s^6 - 16)}{s^2}
					\end{equation*} 
					in $\mathbb{A}^2$. By setting $t = \frac{(s-2)(s+2)(s^2 + 2)^2 u + (s^4 - 2s^2 + 4)}{6s^2}$ we see that $X(G_{11}, \s 3^+)$ is birational to the curve $u^2 = (s^2 - 2s + 4)(s^2 + 2s + 4)$ in $\mathbb{A}^2$. Setting $u = 4(x^2 + 2x + 4)/(x^2 - 4)$ and $s = 4y/(x^2 - 4)$ shows that $X(G_{11}, \s 3^+)$ is isomorphic over $\mathbb{Q}$ to the elliptic curve with Weierstrass equation $y^2 = x^3 + x^2 - 4x - 4$ and LMFDB label \LMFDBLabel{48.a4}.

					The elliptic curve \LMFDBLabel{48.a4} has Mordell-Weil group isomorphic to $(\bbZ/2\bbZ)^2$. From the $j$-maps above it is simple to check that each of these points on $X(G_{11}, \s 3^+)$ is a cusp. In particular the double cover $X(G_{11}(\sqrt{-1}), \s 3)/\delta$ has no noncuspidal rational points.
				\end{proof}
			\end{lemma}
		
			As a corollary we may deduce the level $12$ part of \thref{conj:allqtt}. 
		
			\begin{coro}
				There exist no elliptic curves $E/\bbQ$ admitting a nontrivial quadratic twist-type $(12, 1)$ or $(12, 5)$-congruence.
				
				\begin{proof}
					By \thref{prop:nontrivialCMcongs} no CM elliptic curve $E/\bbQ$ admits a nontrivial quadratic twist-type $12$-congruence. Hence if $E/\bbQ$ admits a nontrivial quadratic twist-type $(12, 1)$ or $(12, 5)$-congruence then by \thref{coro:classify2kpl} the elliptic curve $E$ gives rise to an exceptional point on one of the four modular curves $X(G_9(\sqrt{-1}), \ns 3)/\delta$, $X(G_{11}(\sqrt{-1}), \ns 3)/\delta$, $X(G_9(\sqrt{-1}), \s 3)/\delta$, and $X(G_{11}(\sqrt{-1}), \s 3)/\delta$. But by Lemmas \ref{curves:Gmisc} and \ref{curves:Gmisc4} each of these curves has only rational cusps and CM points. 
				\end{proof}
			\end{coro}

			\begin{lemma} \thlabel{curves:s4ns3}
				The curve $X(\s 4^+, \ns  3^+)$ is isomorphic over $\bbQ$ to the elliptic curve with LMFDB label \LMFDBLabel{36.a4} and rank $0$. Moreover the only rational points on the modular curve $X(\s 4^+, \ns  3^+)$ are cusps and CM points corresponding to the discriminants $-4$ and $-7$. In particular the modular curve $X(\s 4, \ns  3)/\delta$ has no exceptional points.

				\begin{proof}
					To compute a model for $X(\s 4^+, \ns  3^+)$ we fix a model (computed by Rouse and Zureick-Brown \cite{RZB_ECOQA2AIOG}) for $X(\s 4^+) \cong \mathbb{P}^1$ with $j$-invariant map given by $\frac{-64(t^2 - 3)^3 (t^2 + 1)^3}{(t^2 - 1)^4}$. We then see that $X(\s 4^+, \ns  3^+)$ is given by requiring that $(t^2 - 1)$ is a cube. In particular $X(\s 4^+, \ns  3^+)$ is isomorphic over $\bbQ$ the curve 
					\begin{equation*}
						t^2 = s^3 + 1
					\end{equation*} 
					in $\mathbb{A}^2$. This is a Weierstrass equation for the elliptic curve with LMFDB label \LMFDBLabel{36.a4} and rank $0$.

					The elliptic curve \LMFDBLabel{36.a4} has Mordell-Weil group isomorphic to $\bbZ/6\bbZ$. Of these points on $X(\s 4^+, \ns  3^+)$ three are cusps, two correspond to a curve with CM by an order of discriminant $-7$, and one corresponds to a curve with CM by an order of discriminant $-4$. Since $X(\s 4, \ns  3)/\delta$ is a double cover of $X(\s 4^+, \ns  3^+)$ it has no exceptional points.
				\end{proof}
			\end{lemma}

			\begin{lemma} \thlabel{curves:ns4ns5}
				The modular curve $X(\ns 4, \ns 5)/\delta$ has no exceptional points. 

				\begin{proof}
					Baran showed that the genus $2$ curve $X(\ns 4^+, \ns 5^+)$ has no exceptional points (see \cite[Table 5.1]{B_NONSCSMCATCNP}). Hence $X(\ns 4, \ns 5)/\delta$ has no exceptional points.
				\end{proof}
			\end{lemma}

			We now show that the remaining modular curves in Table \ref{table:ratpoints} arise as a double cover of an elliptic curve of rank $1$. The models which we compute are recorded in Table \ref{table:models}. In the cases where we find an exceptional point (i.e., $X(\ns 3, \ns 5)/\delta$ and $X(\ns 8, \ns 3)/\delta$) we also compute a canonical model. Canonical models for the other curves are not recorded here but may be found in the electronic data \cite{ME_ELECTRONIC}.

			\begin{lemma} \thlabel{curves:ns3ns5}
				We have: 
				\begin{enumerate}[label=\eniii]
					\item The modular curve $X(\ns 3^+, \ns 5^+)$ is isomorphic to the elliptic curve with Weierstrass equation $v^2 + v = u^3 + 1$ and LMFDB label \LMFDBLabel{225.c2}.
					\item The modular curve $X(\ns 3, \ns 5)/\delta$ is birational to the curve in $\mathbb{A}^3$ given by 
					\begin{equation} \label{curve:ns3ns5d_affine}
						v^2 + v = u^3 + 1 \qquad \text{and} \qquad w^2 = 3(3u^2 - 2u - 4v + 3)(9u^2 - 10u + 5)
					\end{equation}
					with forgetful map to the model of $X(\ns 3^+, \ns 5^+)$ in \textnormal{(i)} given by $(u,v,w) \mapsto (u,v)$.
					\item The modular curve $X(\ns 3, \ns 5)/\delta$ is a hyperelliptic curve of genus $2$ which is isomorphic over $\bbQ$ to the curve with Weierstrass equation
					\begin{equation} \label{curve:ns3ns5d}
						y^2 + y = x^6 + 3x^5 - 5x^3 + 3x.
					\end{equation}
				\end{enumerate}
	
				\begin{proof}
					(i) We begin by computing equations for $X(\ns 3^+, \ns 5^+)$. Recall the models for $X(\ns 3^+)$ and $X(\ns 5^+)$ with $j$-invariants given by $t^3$ and equation (\ref{eqn:jns5}). It follows that $X(\ns 3^+, \ns 5^+)$ is birational to the curve given by
					\begin{equation*}
						s^3 = \frac{5^3 (t + 1) (2t + 1)^3 (2t^2 - 3t + 3)^3}{(t^2 + t - 1)^5}
					\end{equation*}
					in $\mathbb{A}^2$. Therefore $X(\ns 3^+, \ns 5^+)$ is isomorphic to the elliptic curve with Weierstrass equation $v^2 + v = u^3 + 1$ and LMFDB label \LMFDBLabel{225.c2}. The forgetful maps $X(\ns 3^+, \ns 5^+) \to X(\ns 3^+)$ and $X(\ns 3^+, \ns 5^+) \to X(\ns 5^+)$ are given by 
					\begin{equation*}
						t_3 = \frac{-5 u (v + 3)(8v^2 + 23v + 17) }{(v^2 + v - 1)^2} \qquad \text{and} \qquad t_5 = \frac{-(v + 2)}{(v + 1)}.
					\end{equation*}
					
					(ii) Now $3(t_3^2 + 12t_3 + 144)(8t_5^2 - 12t_5 + 7)$ is equal to $3(3u^2 - 2u - 4v + 3)(9u^2 - 10u + 5)$ up to a square factor in $\bbQ(X(\ns 3^+, \ns 5^+))$. In particular, noting (\ref{eqn:dns3}) and (\ref{eqn:dns5}), $X(\ns 3, \ns 5)/\delta$ is birational to the curve, $\mathcal{C}$, given by (\ref{curve:ns3ns5d_affine}).

					(iii) \texttt{Magma} verifies that taking
					\begin{equation*}
						u = \frac{x^4 + 2x^3 - x + 3y + 4}{(x-1)^2 (x+2)^2}, 
					\end{equation*}
					\begin{equation*}
						v = \frac{x^6 + 3x^5 + 6x^4 + 7x^3 + 3x^2y - 6x^2 + 3xy - 9x + 3y +
						16}{(x-1)^3 (x+2)^3}, 
					\end{equation*}
					\begin{equation*}
						w = \frac{3(2x + 1)(x^6 + 3x^5 + 21x^4 + 37x^3 + 6x^2y - 27x^2 + 6xy - 45x + 15y +
						64)}{(x-1)^4 (x+2)^4}
					\end{equation*}
					gives a birational map between the hyperelliptic curve with Weierstrass equation (\ref{curve:ns3ns5d}) and $\mathcal{C}$.	
				\end{proof}
			\end{lemma}	

			\begin{example}[A $30$-Congruence] \thlabel{ex:30cong}
				The curve $X(\ns {3}, \ns  {5})/\delta$ with model given by (\ref{curve:ns3ns5d}) has the rational point $(3 , 36)$ which lies in the fibre of the non-CM $j$-invariant 
				\begin{equation*}
					\frac{- 3^3 \cdot 5^6 \cdot 199^3 \cdot 809^3 \cdot 5059^3}{61^{15}}.
				\end{equation*}
				The minimal twist with this $j$-invariant is the elliptic curve $E_{(30,7)}/\bbQ$ from the statement of \thref{thm:examples} with conductor $61 \cdot 214663^2$. By \thref{coro:classifyp1p2} the elliptic curve $E_{(30,7)}$ admits a nontrivial $(30,7)$-congruence with its quadratic twist by $-214663$.
			\end{example}
		
			\begin{lemma} \thlabel{curves:ns3s5}
				We have: 
				\begin{enumerate}[label=\eniii]
					\item The modular curve $X(\ns 3^+, \s 5^+)$ is isomorphic to the elliptic curve with Weierstrass equation $v^2 + v = u^3 + 1$ and LMFDB label \LMFDBLabel{225.c2}.
					\item The modular curve $X(\ns 3, \s 5)/\delta$ is birational to the curve in $\mathbb{A}^3$ given by 
					\begin{equation} \label{curve:ns3s5d}
						v^2 + v = u^3 + 1 \qquad \text{and} \qquad w^2 = -(3v + 3u^3 - 12u^2 + 12u - 6)
					\end{equation}
					with forgetful map to the model of $X(\ns 3^+, \s 5^+)$ in \textnormal{(i)} given by $(u,v,w) \mapsto (u,v)$.
				\end{enumerate}	
				
				\begin{proof}
					(i) Recall the model for $X(\s 5^+)$ with $j$-invariant given by equation (\ref{eqn:js5}). Then $X(\ns 3^+, \s 5^+)$ is birational to the curve given by 
					\begin{equation*}
						s^3 = \frac{(t + 3)^3(t^2 + t + 4)^3(t^2- 4t - 1)^3}{(t^2 + t - 1)^5}
					\end{equation*}
					in $\mathbb{A}^2$. It follows that $X(\ns 3^+, \s 5^+)$ is also isomorphic to the elliptic curve with Weierstrass equation $v^2 + v = u^3 + 1$ and LMFDB label \LMFDBLabel{225.c2}. The forgetful maps $X(\ns 3^+, \s 5^+) \to X(\ns 3^+)$ and $X(\ns 3^+, \s 5^+) \to X(\s 5^+)$ are given by 
					\begin{equation*}
						t_3 = \frac{- u (v - 2)(v^2 + v + 4)(v^2 + 6v + 4) }{(v^2 + v - 1)^2} \qquad \text{and} \qquad t_5 = -(v+1).
					\end{equation*}
					
					(ii) Now $-3(t_3^2 + 12t_3 + 144)(t_5^2 - 4t_5 - 16)$ is equal to $-(3v + 3u^3 - 12u^2 + 12u - 6)$ up to a square factor in $\bbQ(X(\ns 3^+, \s 5^+))$. In particular, noting (\ref{eqn:dns3}) and (\ref{eqn:ds5}), $X(\ns 3, \s 5)/\delta$ is birational to the curve given by (\ref{curve:ns3s5d}).
				\end{proof}
			\end{lemma}

		\begin{lemma} \thlabel{curves:ns3ns7}
			We have: 
			\begin{enumerate}[label=\eniii]
				\item The modular curve $X(\ns 3^+, \ns 7^+)$ is isomorphic to the elliptic curve with Weierstrass equation $v^2 + v = u^3 + 12$ and LMFDB label \LMFDBLabel{441.d2}.
				\item The modular curve $X(\ns 3, \ns 7)/\delta$ is birational to the curve in $\mathbb{A}^3$ given by 
				\begin{equation} \label{curve:ns3ns7d}
					v^2 + v = u^3 + 12 \qquad \text{and} \qquad w^2 = 4u^3 + 21u^2 - 42u + 49
				\end{equation}
				with forgetful map to the model of $X(\ns 3^+, \ns 7^+)$ in \textnormal{(i)} given by $(u,v,w) \mapsto (u,v)$.
			\end{enumerate}

			\begin{proof}
				(i) Recall the model for $X(\ns 7^+)$ with $j$-invariant given by equation (\ref{eqn:jns7}). Then $X(\ns 3^+, \ns 7^+)$ is birational to the curve given by 
				\begin{equation*}
					s^3 = \frac{(3t + 1)^3 (t^2 + 3t + 4)^3 (t^2 + 10t + 4)^3 (4t^2 + 5t + 2)^3}{(t^3 + t^2 - 2t - 1)^7}
				\end{equation*}
				in $\mathbb{A}^2$. It follows that $X(\ns 3^+, \ns 7^+)$ is isomorphic to the elliptic curve with Weierstrass equation $v^2 + v = u^3 + 12$ and LMFDB label \LMFDBLabel{441.d2} (see also \cite[Equation 5.5]{B_NONSCSMCATCNP}). The forgetful map $X(\ns 3^+, \ns 7^+) \to X(\ns 7^+)$ is given by
				\begin{equation*}
					t_7 = \frac{6 u^2 - u v - 11 u  - 7 v  + 70 }{3u^2 - u v + 10 u  - 7 v  - 77}.
				\end{equation*}
				Let the forgetful map $X(\ns 3^+, \ns 7^+) \to X(\ns 3^+)$ be given by $(u,v) \mapsto t_3$.
				
				(ii) Now $3(t_3^2 + 12t_3 + 144)(16t_7^4 + 68t_7^3 + 111t_7^2 + 62t_7 + 11)$ is equal to $4u^3 + 21u^2 - 42u + 49$ up to a square factor in $\bbQ(X(\ns 3^+, \ns 7^+))$. In particular, noting (\ref{eqn:dns3}) and (\ref{eqn:dns7}), we see that $X(\ns 3, \ns 7)/\delta$ is birational to the curve given by (\ref{curve:ns3ns7d}).
			\end{proof}
		\end{lemma}
	
		\begin{lemma} \thlabel{curves:2kns11}
			We have: 
			\begin{enumerate}[label=\eniii]
				\item The modular curve $X(\ns  11^+)$ is isomorphic to the elliptic curve with Weierstrass equation $v^2 + v = u^3 - u^2 - 7u + 10$ and LMFDB label \LMFDBLabel{121.b2}.
				\item The modular curve $X(\ns 2, \ns 11)/\delta$ is birational to the curve in $\mathbb{A}^3$ given by 
				\begin{equation} \label{curve:ns2ns11d}
					v^2 + v = u^3 - u^2 - 7u + 10 \qquad \text{and} \qquad w^2 = (-3u^2 + 24u - 37)v - u^4 + 5u^3 + 18u^2 - 95u + 94
				\end{equation}
				with forgetful map to the model of $X(\ns  11^+)$ in \textnormal{(i)} given by $(u,v,w) \mapsto (u,v)$.
			\end{enumerate}
		
			\begin{proof}
				(i) This is  \cite[Theorem~1.1]{CC_ECWNM11R}, \cite{H_SLCMX11}, \cite[Th{\'e}or{\`e}me III.13.18]{H_CESLCM},or \cite[Proposition 4.3.8.1]{L_CMDN11}. The $j$-invariant map is given by 
				\begin{equation*}
					\frac{(f_1 f_2 f_3 f_4)^3}{f_5^2 f_6^{11}}
				\end{equation*}
				where 
				\begin{equation*}
					\begin{aligned}
						&f_1 = u^2 + 3u + 6 ,\\
						&f_2 = 11(u^2 - 5)v + (2u^4 + 23u^3 - 72u^2 - 28u + 127) ,\\
						&f_3 = 6v + 11u - 19,\\
						&f_4 = 22(u - 2)v + (5u^3 + 17u^2 - 112u + 120),\\
						&f_5 = 11v + (2u^2 + 17u - 34), \\
						&f_6 = (u - 4)v - (5u - 9).
					\end{aligned}
				\end{equation*}
				
				(ii) Let $\mathcal{E}/\bbQ(X(\ns  11^+))$ be an elliptic curve with this $j$-invariant. By \thref{rem:howtogetd} and \cite[Proposition III.13.25]{H_CESLCM} we see that $\mathcal{E}$ is $(11, 1)$-congruent to its quadratic twist by $d = d_1 d_2 d_3$ where 
				\begin{equation*}
					\begin{aligned}
						&d_1 = 3619v - (108u^5 + 168u^4 - 641u^3 - 1007u^2 - 3769u + 9670), \\
						&d_2 = 48(3u + 5)v - (276u^2 - 452u + 185), \\
						&d_3 = 16v + 28u - 49.
					\end{aligned}
				\end{equation*}
				
				Now $d \Delta(\mathcal{E})$ is equal to $(-3u^2 + 24u - 37)v - u^4 + 5u^3 + 18u^2 - 95u + 94$ up to a square factor in $\bbQ(X(\ns 11^+))$ so $X(\ns  2, \ns  11)/\delta$ is birational to the curve given by (\ref{curve:ns2ns11d}). 
			\end{proof}
		\end{lemma}

		\begin{lemma} \thlabel{curves:ns8ns3}
			We have: 
			\begin{enumerate}[label=\eniii]
				\item The modular curve $X(\ns 8^+, \ns 3^+)$ is isomorphic to the elliptic curve with Weierstrass equation $v^2 = u^3 + 8$ and LMFDB label \LMFDBLabel{576.e4}.
				\item The modular curve $X(\ns 8, \ns 3)/\delta$ is birational to the curve in $\mathbb{A}^3$ given by 
				\begin{equation} \label{curve:ns8ns3d_affine}
					v^2 = u^3 + 8 \qquad \text{and} \qquad w^2 = 3(3u^3 + 3u^2 + 2uv - 8u + 8)
				\end{equation}
				with forgetful map to the model of $X(\ns 8^+, \ns 3^+)$ in \textnormal{(i)} given by $(u,v,w) \mapsto (u,v)$.
				\item The modular curve $X(\ns  8, \ns  3 )/\delta$ is a non-hyperelliptic curve of genus $3$ which may be canonically embedded in $\bbP^2$ over $\bbQ$ as the curve
				\begin{equation} \label{eqn:48curve}
					\begin{aligned}
						2X^3Y - 2X^3Z + X^2Y^2 - 2X^2YZ - 2X^2Z^2 - 2XY^3 + 2XY^2Z \\+ 2XYZ^2 - 2XZ^3 - Y^4 + 2Y^3Z + Y^2Z^2 - 2YZ^3 = 0.
					\end{aligned}
				\end{equation}
			\end{enumerate}

			\begin{proof}
				(i) We fix a model for $X(\ns  8^+) \cong \bbP^1$ so that the $j$-invariant map is given by 
				\begin{equation*}
					\frac{-2^{15}(t -2)t^3(t^2 - 2t - 2)^3(t^2 - 2t + 2)^3}{(t^2 - 2)^8},
				\end{equation*}
				this is the model from \cite{RZB_ECOQA2AIOG} obtained after replacing $t$ by $1/t$.
				
				Then $X(\ns  8^+, \ns  3^+ )$ is birational to the curve given by 
				\begin{equation*}
					s^3 = \frac{-2^{15}(t -2)t^3(t^2 - 2t - 2)^3(t^2 - 2t + 2)^3}{(t^2 - 2)^8}
				\end{equation*}
				in $\mathbb{A}^2$. It follows that $X(\ns 8^+, \ns 3^+)$ is isomorphic to the elliptic curve with Weierstrass equation $v^2 = u^3 + 8$ and LMFDB label \LMFDBLabel{576.e4}. The forgetful maps $X(\ns 8^+, \ns 3^+) \to X(\ns 8^+)$ and $X(\ns 8^+, \ns 3^+) \to X(\ns 3^+)$ are given by
				\begin{equation*}
					t_8 = \frac{2(v - 2)}{v - 4} \qquad \text{and} \qquad t_3 = \frac{64 u (v - 2) (v^2 - 12v + 24) (v^2 - 4v + 8)}{(v^2 - 8)^3}.
				\end{equation*}
				
				(ii) Let $\mathcal{E}/\bbQ(X(\ns 8^+, \ns 3^+))$ be an elliptic curve with $j$-invariant $t_3^3$. Then $-3(t_3^2 + 12t_3 + 144)\Delta(\mathcal{E})$ is equal to $3(3u^3 + 3u^2 + 2uv - 8u + 8)$ up to a square factor in $\bbQ(X(\ns 8^+, \ns 3^+))$. Therefore, recalling (\ref{eqn:dns3}), $X(\ns 8, \ns 3)/\delta$ is birational to the curve $\mathcal{C}$ given by (\ref{curve:ns8ns3d_affine}).
				
				(iii) \texttt{Magma} verifies that taking
				\begin{equation*}
					u = \frac{2(XY + Y^2 - XZ - YZ - Z^2)}{Z^2}, 
				\end{equation*}
				\begin{equation*}
					v = \frac{2(2X^2Y + 2XY^2 - 2X^2Z - 3XYZ - Y^2Z - 2XZ^2 + YZ^2)}{Z^3}, 
				\end{equation*}
				\begin{equation*}
					w = \frac{2(X + 2Y - Z)(2XY + 2Y^2 - 2XZ - 2YZ - 3Z^2)}{Z^3}
				\end{equation*}
				gives a birational map between the curve given by (\ref{eqn:48curve}) and $\mathcal{C}$.
			\end{proof}
		\end{lemma}
		
		\begin{example}[A $48$-Congruence] \thlabel{ex:48cong}
			The model for $X(\ns  {8}, \ns  {3})/\delta$ given by (\ref{eqn:48curve}) has the rational point $(4 : -7 : 1)$ which lies in the fibre of the non-CM $j$-invariant 
			\begin{equation*}
				\frac{2^{15} \cdot 3^3 \cdot 5^3 \cdot 31^3 \cdot 41^3 \cdot 47^3 \cdot 83^3 \cdot 293^3}{23^{24}}.
			\end{equation*}
			The minimal twist with this $j$-invariant is the elliptic curve $E_{(48,19)}/\bbQ$ from the statement of \thref{thm:examples}. The conductor of $E_{(48,19)}$ is $19^2 \cdot 23 \cdot 43^2 \cdot 39451^2$ and its discriminant is $-19^3 \cdot 23^{24} \cdot 43^3 \cdot 39451^3$. 
			
			By \thref{coro:classify2kpl} the elliptic curve $E_{(48,19)}$ admits a nontrivial $(48,19)$-congruence with its quadratic twist by its discriminant, which (up to a square factor) is equal to $-19\cdot 43 \cdot 39451$.
		\end{example}

		\subsection{A Review of the Mordell-Weil Sieve} \label{sec:MWsieve}
			In order to search for points on the modular curves in Table \ref{table:models} and prove \thref{thm:conditional} we briefly recall the Mordell-Weil sieve (see also \cite[Section 10]{RSZB_LAIOGFECOQ}). 
			
			Let $\phi \colon C \to E$ be a morphism defined over $\bbQ$ where $E/\bbQ$ is an elliptic curve and $C/\bbQ$ has genus $\geq 2$. Let $P_1, ..., P_n \in E(\bbQ)$ be generators for the Mordell-Weil group of $E$. For any good prime $p$ for both $C$ and $E$ we have a commutative diagram
			\[\begin{tikzcd}
				{C(\bbQ)} & {E(\bbQ)} & {\bbZ^n} \\
				{\widetilde{C}(\bbF_p)} & {\widetilde{E}(\bbF_p)} & {\prod_{i=1}^{n} \bbZ/N_{p,i}\bbZ}
				\arrow["\operatorname{red}_p", from=1-2, to=2-2]
				\arrow[two heads, from=1-3, to=2-3]
				\arrow[from=2-3, to=2-2]
				\arrow[two heads, from=1-3, to=1-2]
				\arrow["\phi", from=1-1, to=1-2]
				\arrow[from=1-1, to=2-1]
				\arrow["{\widetilde{\phi}}", from=2-1, to=2-2]
			\end{tikzcd}\]
			where $N_{p, i}$ is the order of $\operatorname{red}_p(P_i)$ in $\widetilde{E}(\bbF_p)$, and $\widetilde{C}$, $\widetilde{E}$, $\widetilde{\phi}$ are the reductions of $C$, $E$, and $\phi$ modulo $p$ respectively. If there exists a rational point in the fibre of $\phi$ above $Q \in E(\bbQ)$ then the fibre of $\widetilde{\phi}$ above $\operatorname{red}_p(Q)$ is non-empty. Writing $Q = \sum_{i=1}^{n} a_i P_i$ we obtain a set of congruence conditions on $a = (a_1,...,a_n)$ in $\prod_{i=1}^{n} \bbZ/N_{p,i}\bbZ$.

			Repeating this for a (finite) set of primes of good reduction for $E$ and $C$ allows us to give an (often large) bound, $B > 0$, such that every rational point on $C$ is either known to us, or lies in the fibre above a point $Q = \sum_{i=1}^{n} a_i P_i$ such that $\max \{|a_i|\} _{i=1}^{n} > B$.
			
			We are now able to prove \thref{thm:conditional}.
			
			\begin{proof}[Proof of \thref{thm:conditional}]
				Suppose that $E/\bbQ$ admits a nontrivial quadratic twist-type $(N,r)$-congruence. If $N$ is a prime power, assumption (A) together with \thref{thm:infmanypl} shows that either $(N, r) \in \mathcal{S}$ or $E$ is a quadratic twist of either $E_{(32,3)a}$ or $E_{(32,3)b}$. 
				
				It therefore suffices to consider the case where $N$ is divisible by two or more distinct primes. If $p$ is odd, the case when $N = 2p^l$ follows immediately. Otherwise by assumptions (A) and (C), together with Propositions \ref{prop:nontrivialCMcongs} and \ref{thm:infmanypl}, and Corollaries \ref{coro:classifyp1p2} and \ref{coro:classify2kpl}, if $(N,r) \not\in \mathcal{S}$ then there exists an exceptional point on one of the $13$ curves in Table \ref{table:ratpoints} above $j(E)$.

				We proved in Lemmas \ref{curves:2ks3}, \ref{curves:2ks5}, \ref{curves:Gmisc}, \ref{curves:Gmisc4}, \ref{curves:s4ns3}, and \ref{curves:ns4ns5} that the modular curves appearing in Table \ref{table:ratpoints} but not Table \ref{table:models} have no exceptional points. Hence $j(E)$ must be an exceptional $j$-invariant for one of the $5$ curves in Table \ref{table:models}.
				
				The Mordell-Weil sieve shows that, except for the points on $X(\ns 3, \ns 5)/\delta$ and $X(\ns 8, \ns 3)/\delta$ above $j(E_{(30,7)})$ and $j(E_{(48,19)})$, the only exceptional points on these $5$ curves lie in the fibre of $X(H_1, H_2)/\delta \to X(H_1^+, H_2^+)$ above a point of canonical height $\geq 10^{15}$.
			\end{proof}
				
	\section{Quadratic Twist-Type \texorpdfstring{$(15, 7)$}{(15,7)}-Congruences and Rational Points on \texorpdfstring{${X(\ns 3, \ns 5)/\delta}$}{X(ns 3, ns 5)/delta}} \label{sec:15-7ratpts}

		Our final aim is to study the rational points on the modular curve ${X(\ns 3, \ns 5)/\delta}$. In particular we will show that (up to a quadratic twist) $E_{(30,7)}$ is the only elliptic curve over $\bbQ$ admitting a nontrivial quadratic twist-type $(15,7)$-congruence. The key observation is that since the Jacobian of ${X(\ns 3, \ns 5)/\delta}$ admits a Richelot isogeny we may use elliptic Chabauty to determine the rational points on ${X(\ns 3, \ns 5)/\delta}$.   
			
		More specifically we will prove: 			
						
		\begin{prop} \thlabel{prop:ns3ns5pts}
			The rational points on the model for $X(\ns 3, \ns 5)/\delta$ given by (\ref{curve:ns3ns5d}) are exactly those listed in Table \ref{table:pts_ns3ns5d}. In particular the only nontrivial quadratic twist-type $(15, 7)$-congruence over $\bbQ$ (up to taking quadratic twists) is between the elliptic curve $E_{(30,7)}$ (see \thref{thm:examples}) and its quadratic twist by $-214663 $.
		\end{prop}
		
		\begin{table}[H]
			\centering
			\begin{tabular}{c|c|c|c}
				$j$-Invariant    & Points       & Quadratic Twist       & CM     \\
				\hline
				$0$ & $ (1 , 1), (-2 , 1) $ & $-3 $ & $-3$ \\
				$-3^{3} \cdot 5^{3} $ & $ (-1 , 0 ), (0 , 0 ) $ & $-7 $ & $-7$ \\
				$3^{3} \cdot 5^{3} \cdot 17^{3} $ & $ (1 , -2 ), (-2 , -2 ) $ & $-7 $
				& $-28$ \\
				$-2^{18} \cdot 3^{3} \cdot 5^{3} $ & $+\infty, -\infty$ & $-43 $
				& $-43$ \\
				$-2^{15} \cdot 3^{3} \cdot 5^{3} \cdot 11^{3} $ & $(-1 , -1 ), (0 , -1 )$ & $-67 $ & $-67$ \\
				$-2^{18} \cdot 3^{3} \cdot 5^{3} \cdot 23^{3} \cdot 29^{3} $ & $ (3 , -37 ), (-4 , -37 ) $ & $-163 $ & $-163$ \\
				$-\frac{3^{3} \cdot 5^{6} \cdot 199^{3} \cdot 809^{3} \cdot 5059^{3} }{61^{15}
				}$ & $ (3 , 36 ), (-4 , 36 ) $ & $-214663 $ &  \\
			\end{tabular}
			\caption{All rational points on the model for $X(\ns 3, \ns 5)/\delta$ given by (\ref{curve:ns3ns5d}).} \label{table:pts_ns3ns5d}
		\end{table}
			
		\subsection{\texorpdfstring{$2$}{2}-Cover Descent on Genus \texorpdfstring{$2$}{2} Curves Admitting a Richelot Isogeny} \label{sec:ratpts}	

			Let $C/\bbQ$ be a hyperelliptic curve of genus $2$ given by a Weierstrass equation $y^2 = f(x)$, where $f(x) \in \bbQ[x]$ has degree $6$. Further suppose that there exists an extension $K/\bbQ$ such that $f(x)$ factors as a product $f_1(x) f_2(x)$ over $K$ where $f_1(x)$ has degree $4$ and $f_2(x)$ has degree $2$. 

			We consider the family of genus $3$ double covers $C_d \to C$ over $K$, where $C_d$ is given by 
			\begin{align*}
				d y_1^2 &= f_1(x) \\
				d y_2^2 &= f_2(x).
			\end{align*}
			for some $d \in K^\times$. Let $S$ be the union of the infinite places of $K$, the primes of $K$ lying above $2$ and the bad primes of $C$. We write
			$$K(S,2) = \{a \in K^\times /(K^\times)^2 : \text{$v_{\mathfrak{p}}(a) \equiv 0 \pmod{2}$ for each (finite) place $\mathfrak{p}$ of $K$ not in $S$} \}$$ 
			where $v_{\mathfrak{p}} \colon K_{\mathfrak{p}} \to \mathbb{Z}$ is the normalised discrete valuation.

			By {\'e}tale descent (see e.g., \cite[Lemma 3.1]{B_CMUEC}), for each point $P \in C(K)$ there exists a $d \in K(S,2)$ such that $P$ lifts to a $K$-point on $C_d$. However, a rational point $P \in C(\bbQ)$ lifts to a point on $C_d(K)$ with rational $x$-coordinate. We therefore consider the family of genus $1$ curves $C'_d / K$ given by 
			\begin{equation*}
				C'_d : d y^2 = f_1(x).
			\end{equation*}

			Thus, to determine the $x$-coordinates supporting a rational point on $C$ it suffices to determine the set of $P \in C'_d(K)$ such that $x(P) \in \bbQ$. If $C_d'$ has a $K$-point then the method of elliptic Chabauty provides a solution to this problem under the hypothesis that $\operatorname{rank} C_d'(K) \leq [K : \bbQ]$ (this has been implemented in \texttt{Magma}, see \cite{B_CMUEC} and \cite{B_STDEOSnn2}).

			Now $\operatorname{Jac}(C)$ admits a Richelot isogeny over $\bbQ$ if and only if there exists a degree $3$, {\'e}tale, $\bbQ$-algebra $L/\bbQ$ such that $f(x)$ admits such a factorisation $f(x) = \operatorname{Norm}_{L[x]/\bbQ[x]} Q(x)$ for some quadratic polynomial $Q(x) \in L[x]$ (see \cite[Lemma 4.1]{BD_TAOG2CW44SJ}). In particular, for each number field $L'$ appearing in the factorisation of $L$ we obtain a factorisation $f(x) = f_1(x) f_2(x)$ defined over $L'$ (but not over any proper subfields). 
			\begin{remark}
				In the file \texttt{twocover.m} in the supporting electronic data \cite{ME_ELECTRONIC} we have implemented the above method for genus $2$ curves $C/\mathbb{Q}$ given by a Weierstrass equation $y^2 = f(x)$ where $f(x)$ is irreducible and $\operatorname{Jac}(C)$ admits a Richelot isogeny.
			\end{remark}		
		
			We now specialise to the case of rational points on $X(\ns 3, \ns 5)/\delta$ as follows. Let $C$ be the genus $2$ curve with short Weierstrass equation $$y^2 = f(x) = 4x^6 + 12x^5 - 20x^3 + 12x + 1,$$ and note that replacing $y$ by $2y + 1$ gives an isomorphism from the model for $X(\ns 3, \ns 5)/\delta$ given by (\ref{curve:ns3ns5d}) to $C$. The \texttt{Magma} function \texttt{RankBounds} verifies that $J = \operatorname{Jac}(C)$ has rank $2$. In particular Chabauty's method is not immediately applicable since $\operatorname{rank} J(\bbQ) = \dim J$.  
			
			Let $\theta$ be a root of $x^3 + x^2 - 3x + 3$ and let $K$ be the cubic number field $\bbQ(\theta)$. Over $K$, $f(x)$ factors as a product of a quartic and quadratic polynomial $f_1(x)$ and $f_2(x)$, where
			\begin{align*}
				f_1(x) &= 4 x^4 + 8 x^3 + (-2\theta^2 - 2\theta) x^2 + (-2\theta^2 - 2\theta - 4) x + \theta^2 - 2\theta + 1 \\
				f_2(x) &=  x^2 + x + (\theta^2 + \theta - 4)/2.
			\end{align*}

			Let $S$ be the union of the set of infinite places of $K$ and set of places of $K$ lying above $2, 3,$ or $5$. We consider the family of $64$ genus $3$ double covers $C_d \to C$ over $K$ indexed by $d \in K(S,2)$, and the corresponding family of genus $1$ curves $C_d'$. 
			
			Of these $64$ genus $1$ curves, $60$ have a local obstruction at a place $v \in S$ - that is $C_d'(K_v) = \emptyset$ for some $v \in S$. The remaining four curves are when $d \in \{ 1, \theta^2 + 2\theta - 1, \theta^2 + \theta - 3, 2\theta^2 + 2\theta - 3 \}$. Each of these four curves $C_d'$ have a $K$-rational point and have Mordell-Weil rank $\leq 2$ - in particular elliptic Chabauty may be applied (implemented in \texttt{Magma} as the function \texttt{Chabauty}). We find that the $x$-coordinates supporting a rational point on $C$ are contained in the set $\{-4, -2, -1, -1/2, 0, 1, 3, \infty\}$. Each of these $x$-coordinates, except $-1/2$, gives rise to a rational point on $C$. \thref{prop:ns3ns5pts} then follows from an elementary calculation.

		\subsection{A Solution to the Class Number \texorpdfstring{$1$}{1} Problem} \label{sec:clnum}

			It was conjectured by Gauss that the only orders in an imaginary quadratic field with class number $1$ are those with discriminants $-3$, $-4$, $-7$, $-8$, $-11$, $-12$, $-16$, $-19$, $-27$, $-28$, $-43$, $-67$, and $-163$. It is well know that there is a bijective correspondence between orders, $\mathcal{O}$, in an imaginary quadratic field with class number $1$ and $j$-invariants of elliptic curves $E/\bbQ$ with CM by $\mathcal{O}$ (see \cite[Chapter 2]{S_ATITAOEC}). In this case $j(E)$ is an integer (see \cite[Chapter 2 Section 6]{S_ATITAOEC}). 
			
			Since CM elliptic curves give rise to rational points on the modular curve $X(\ns {p^{l}}^+)$ for suitable $p^l$ (see \cite[Section 4]{B_NONSCSMCATCNP} and \cite[A.5]{S_LOTMWT}) this suggests that to prove the class number $1$ problem it is enough to find all the integral points (that is, points lying above $j$-invariants which are integers) on $X(\ns {p_1^{l_1}}^+, ..., \ns {p_n^{l_n}}^+)$ where $p_1,..., p_n$ are distinct prime numbers. This approach has successfully been used to give a solutions to the class number one problem in all cases where the genus of $X(\ns {p_1^{l_1}}^+, ..., \ns {p_n^{l_n}}^+)$ is $\leq 2$ (see \cite{B_NONSCSMCATCNP}, \cite[p. 197]{S_LOTMWT}, and the references therein). 

			We now discuss a slightly different (but closely related) approach to the class number $1$ problem via the modular curves $X(\ns {p^{l}}, \ns {q^{k}})/\delta$ constructed in Section \ref{sec:doublecovers}. We are careful not to quote the results from Section \ref{sec:trivialCM} since they rely upon the classification of elliptic curves $E/\bbQ$ with CM.

			Suppose that an elliptic curve $E/\bbQ$ with $j(E) \neq 0, 1728$ has CM by an order in an imaginary quadratic field $K$. Let $-D$ be the discriminant of $K$, where $D$ is a positive integer. Let $p$ be a prime number and let $l \geq 1$ be an integer. If $p$ is inert in $K$ then $\bar{\rho}_{E, p^l}(G_{\bbQ})$ is contained in the normaliser of a nonsplit Cartan subgroup and $\bar{\rho}_{E, p^l}(G_{\bbQ(\sqrt{-D})})$ is contained in the corresponding Cartan subgroup (see \cite[Proposition 4.1]{B_NONSCSMCATCNP} and \cite[A.5]{S_LOTMWT}). By \thref{cartan_cong_pl} $E$ is therefore $(p^l, -\xi)$-congruent to its quadratic twist by $-D$.
			
			Each prime number $p < \frac{1 + D}{4}$ is inert in $K$ (see the discussion following Proposition 4.10 in \cite{B_NONSCSMCATCNP} and \cite[p. 190]{S_LOTMWT}). Hence if an elliptic curve $E/\bbQ$ has CM by an order in an imaginary quadratic field $K$ of discriminant $-D \leq -4p$, then $E$ is $(p^l, -\xi)$-congruent to its quadratic twist by $-D$.
	
			Now fix a pair of prime powers $p^l, q^k$. By the preceding discussion if $E/\bbQ$ has CM by an order in an imaginary quadratic field $K$ of discriminant $-D \leq \min\{-4p, -4q\}$ then there exists a integral point in the fibre of the $j$-map $X(\ns {p^{l}}^+, \ns {q^{k}}^+) \to X(1)$ above $j(E)$. 
			
			But as we observed, $E$ is $p^{l} q^{k}$-congruent to $E^{-D}$. Therefore by \thref{prop:Xcong} there exists a integral point on $X(\ns {p^{l}}, \ns {q^{k}})/\delta$ above $j(E)$.
			
			In \thref{prop:ns3ns5pts} we determined all the rational (hence integral) points on $X(\ns 3, \ns 5)/\delta$. By taking $p^{l} = 3$ and $q^{k} = 5$ in the above discussion it follows from \thref{prop:ns3ns5pts} that the only imaginary quadratic fields with class number $1$ and discriminant $\leq -20$ have discriminant $-43$, $-67,$ and $-167$. Hence we have a solution to the class number one problem.
		
			\begin{remark}
				Level $15$ proofs of the class number $1$ problem have previously been given by Siegel \cite{S_ZBDSS} and I. Chen \cite{C_OSMCOL15ATCNOP} by analysing integral points on the rank $1$ elliptic curve $X(\ns 3^+, \ns 5^+)$. We have instead analysed the rational points on the double cover $X(\ns 3, \ns 5)/\delta$.
			\end{remark}
	
			It would be interesting to give proofs of the class number $1$ problem by analysing the integral points on the modular curves $X(\ns 2^k, \ns p^l)/\delta$ of genus $0$ and $1$ computed in Section \ref{sec:2kpl}.

\appendix	\section{Tables}

	By Propositions \ref{cartan_cong_pl}, \ref{prop:2adiccongs}, and \ref{thm:infmanypl} we have a finite list of pairs of subgroups $H, H^+ \subset \GL_2(\bbZ/p^l\bbZ)$ which give rise to nontrivial quadratic twist-type $p^l$-congruences where $p^l$ is a prime power (assuming there are no exceptional points on $X(\ns {p^l}^+)$ for $p^l \geq 19$). 
	
	We therefore compute the genera of the modular curves $X(H_1, H_2)/\delta$ which give rise to nontrivial quadratic twist-type $p_1^{l_1}p_2^{l_2}$-congruences over $\bbQ$. In Table \ref{table:pqgenera} we list these genera when $p_1$ and $p_2$ are odd and in Table \ref{table:2kgenera} we list these genera in the case where $p_1 = 2$. To perform these calculations we use Rouse, Sutherland, and Zureick-Brown's function \texttt{GL2Genus} in the package \texttt{gl2.m} \cite{RSZB_Electronic}.

	In Table \ref{table:examples} we give an example of a quadratic twist-type $(N, r)$-congruence for each pair $(N, r)$ allowed by \thref{conj:allqtt}.

	\begin{table}[ht]
		\centering
		\begin{tabular}{c|ccccc}
			\diagbox{$H_1$}{$H_2$}   	     	& $\ns 5$		& $\s 5$ 		& $\ns 7$		& $\s 7$		& $\ns 11$		\\ 
			\hline
			$\ns 3$  							& $\mathit{2}$	& $\mathit{3}$	& $\mathit{4}$	& $7$			& $17$			\\
			$\s 3$   							& $5$			& $8$			& $12$			& $17$			& $40$			\\
			$\ns 5$  							&				&				& $26$			& $35$			& $76$			\\
			$\s 5$   							&				&				& $40$			& $55$			& $117$			\\
			$\ns 7$  							& 				&				& 				&				& $166$			\\
			$\s 7$   							&				&				&				&				& $225$ 	
			
		\end{tabular}

		\caption{Genera of the modular curves $X(H_1, H_2)/\delta$ in \thref{coro:classifyp1p2} of level $p q$ such that $X(H_i^+)$ has infinitely many rational points for $i=1,2$. Models for the curves whose genera are in \textit{italics} are computed in Section \ref{sec:bigconj}.} \label{table:pqgenera}
	\end{table}

	\begin{table}[ht]
		\centering
		\begin{tabular}{c|cccccccc}
		\diagbox{$H_1$}{$H_2$} 		& $\ns 3$	 	& $\s 3$ 		& $\ns 5$ 		& $\s 5$ 		& $\ns 7$ 		& $\s 7$ 		& $\ns 9$ 		& $\ns 11$  \\ 
		\hline
		$G_{9}(\sqrt{-1})$			& $\mathbf{1}$	& $\mathbf{1}$	& $5$			& $7$			& $13$			& $17$			& $19$			& $41$        \\
		$G_{11}(\sqrt{-1})$			& $\mathbf{1}$	& $\mathit{3}$	& $7$			& $9$			& $15$			& $21$			& $21$			& $45$        \\
		\hdashline
		$G(\ns 2)$					& $\mathbf{0}$	& $\mathbf{0}$	& $\mathbf{0}$	& $\mathbf{1}$	& $\mathbf{1}$	& $\mathbf{1}$	& $\mathbf{1}$	& $\mathit{4}$\\
		$G_{9}(\sqrt{-\Delta})$		& $2$			& $3$			& $7$			& $10$			& $16$			& $21$			& $22$			& $46$        \\
		$G_{11}(\sqrt{-\Delta})$	& $0$			& $1$			& $5$			& $8$			& $12$			& $17$			& $16$			& $40$        \\
		\hline
		$G_{97}(\sqrt{-1})$ 		& $7$ 			& 				&				&				&				&				& 				&             \\
		\hdashline
		$G(\ns 4)$					& $\mathbf{0}$	& $\mathbf{1}$	& $\mathit{3}$	& $6$			& $8$			& $11$			& $10$			&             \\
		$G_{90}(\sqrt{-\Delta})$	& $5$			& $13$			&				& 				&				& 				&				&             \\
		\hdashline
		$G_{63}(\sqrt{-1})$			& $7$			&				& 				& 				&  				& 				& 				&             \\
		\hdashline
		$G(\s 4)$					& $\mathit{2}$	& $5$ 			& $13$			& $20$			& $30$			& $41$ 			& $40$			&             \\
		$G_{91}(\sqrt{-\Delta})$	& $5$			& $13$			& 				& 				& 				& 				&				&             \\

		\hline
		$G(\ns 8)$					& $\mathit{3}$	& $9$			&				&				&				& 				& 				&            
		\end{tabular}

		\caption{Genera of the modular curves $X(H_1, H_2)/\delta$ in \thref{coro:classify2kpl} of level $2^k p^l$ such that $X(H_i^+)$ has infinitely many rational points for $i=1,2$. If the image of $X(H_1, H_2)/\delta$ at level $2^{k-1} p^l$ is of genus $\geq 2$ we leave a space. Models for the curves whose genera are in \textbf{bold} are computed in Section \ref{sec:2kpl} and models for the curves whose genera are in \textit{italics} are computed in Section \ref{sec:bigconj}.} \label{table:2kgenera}
	\end{table}

	\begin{landscape} 
	\hspace{0pt}
	\vfill

	\begin{table}[ht]
		\centering
		\begin{tabular}{c|ccccc}
			$(N, r)$	& {\makecell[c]{$j$-invariant of an elliptic curve $E/\bbQ$ admitting \\a quadratic twist-type $(N,r)$-congruence}} 		
																													& {\makecell[c]{LMFDB label of the  \\minimal twist of $E$}} 				
																																					& Quadratic twist	& Modular Curve				& Model\\
			\hline
			$(8, 1)$	& $\frac{-2^6 \cdot 3^3 \cdot 17^3 \cdot 47^3 }{7^8}$										& \LMFDBLabel{16128.w2}				& $-1$				& $X(G_{97})$	& \cite[\href{https://users.wfu.edu/rouseja/2adic/X97.html}{$X_{97}$}]{RZB_Electronic}\\
			$(8, 5)$	& $\frac{2^6 \cdot 17^3 \cdot 19^3 }{3^4 }$													& \LMFDBLabel{2400.p1}				& $-1$				& $X(G_{63})$ 	& \cite[\href{https://users.wfu.edu/rouseja/2adic/X63.html}{$X_{63}$}]{RZB_Electronic}\\
			$(8, 7)$	& $\frac{2^2 \cdot 3^3 \cdot 13^3 }{5^4 }$													& \LMFDBLabel{40.a4}					& $\Delta(E)$		& $X(\s 4^+)$	& \cite[\href{https://users.wfu.edu/rouseja/2adic/X23.html}{$X_{23}$}]{RZB_Electronic}\\
			$(10, 1)$	& $\frac{-2^{9} \cdot 3^3 \cdot 7^3 }{19^5 }$												& \LMFDBLabel{608.d1}				& $-1$				& $X(\s 5^+)$				& \cite[Th{\'e}or{\`e}me III.5.9]{H_CESLCM}\\
			$(12, 11)$	& $\frac{5^3 \cdot 31^3 }{2^6 \cdot 3^3}$													& \LMFDBLabel{726.b2}				& $\Delta(E)$		& $X(\ns 4, \s 3)/\delta$	& \thref{curves:2ks3} \\
			$(20, 3)$	& $\frac{-2^{18} \cdot 3 \cdot 5^3 \cdot 11^3 \cdot 17 \cdot 23^3 \cdot 31^3}{61^{10}}$		& \LMFDBLabel{475983.f1}			& $\Delta(E)$		& $X(\ns 4, \ns 5)/\delta$	& \thref{curves:2kns5} \\
			$(22, 1)$	& $\frac{-2^{9} \cdot 3^3 \cdot 5^3 \cdot 13 \cdot 71^3 \cdot 181^3}{43^{11}}$				& \LMFDBLabel{232544.i1}		& $-1$				& $X(\ns 11^+)$				& \cite{H_SLCMX11} \\
			$(28, 3)$	& $\frac{2^{18} \cdot 3^3 \cdot 5^3 \cdot 7^4 \cdot 37^3 }{13^{14}}$						& \texttt{24721333*}			& $\Delta(E)$		& $X(\ns 2, \ns 7)/\delta$ 	& \thref{curves:2k7}\\
			$(28, 11)$	& $\frac{-2^{15} \cdot 5^3 \cdot 11^3 \cdot 17^3 \cdot 29^3}{13^{14}}$						& \texttt{1218886357*}			& $\Delta(E)$		& $X(\ns 2, \s 7)/\delta$ 	& \thref{curves:2k7}\\
			$(30, 7)$	& $\frac{- 3^3 \cdot 5^6 \cdot 199^3 \cdot 809^3 \cdot 5059^3}{61^{15}}$					& \texttt{2810892417709*}		& $-214663$			& $X(\ns 3, \ns 5)/\delta$	& \thref{curves:ns3ns5}\\
			$(32, 3)$	& $\frac{-2^{18} \cdot 3 \cdot 5^3 \cdot 13^3 \cdot 41^3 \cdot 107^3}{17^{16}}$				& \texttt{356257899*}			& $\Delta(E)$		& $X(\ns 16^+)$ 	& \cite[\href{https://users.wfu.edu/rouseja/2adic/X441.html}{$X_{441}$}]{RZB_Electronic}\\
			$(36, 19)$	& $\frac{-2^{24} \cdot 3^3 \cdot 5^3 \cdot 11^3 \cdot 47^3 \cdot 103^3}{17^{18}}$			& \texttt{21043447153673*}		& $\Delta(E)$		& $X(\ns 2, \ns 9)/\delta$	& \thref{curves:2kns9}\\
			$(48, 19)$	& $\frac{2^{15} \cdot 3^3 \cdot 5^3 \cdot 31^3 \cdot 41^3 \cdot 47^3 \cdot 83^3 \cdot 293^3}{23^{24}}$
																													& \texttt{23893951694358047*}	& $\Delta(E)$		& $X(\ns 8, \ns 3)/\delta$	& \thref{curves:ns8ns3}
		\end{tabular}

		\caption{All known pairs $(N, r)$ for which we have examples of nontrivial quadratic twist-type $(N, r)$-congruences. If $M$ divides $N$ and $(N, r)$ appears in the above list, we do not give an example of an $(M, r)$-congruence. We also list the modular curve giving rise to the corresponding example. When the conductor of an elliptic curve is outside the range of the LMFDB we use a ``$*$'' to replace the label for the isogeny class.} \label{table:examples}
	\end{table}
	\vfill
	\end{landscape}

	\newpage
\newcommand{\etalchar}[1]{$^{#1}$}
\providecommand{\bysame}{\leavevmode\hbox to3em{\hrulefill}\thinspace}
\providecommand{\MR}{\relax\ifhmode\unskip\space\fi MR }
\providecommand{\MRhref}[2]{%
	\href{http://www.ams.org/mathscinet-getitem?mr=#1}{#2}
}
\providecommand{\href}[2]{#2}

\end{document}